\documentclass[final]{siamltex}
\usepackage{epsfig}
\usepackage{graphicx}
\usepackage{amsmath}
\usepackage{amssymb}
\usepackage{setspace}
\usepackage{graphicx}
\usepackage{graphicx}

\title{Landmark and Intensity Based Registration with Large Deformations via Quasi-conformal Maps}

\author{Ka Chun Lam and Lok MIng Lui}

\begin{document}

\maketitle

\begin{abstract}
Registration, which aims to find an optimal one-to-one correspondence between different data, is an important problem in various fields. This problem is especially challenging when large deformations occur. In this paper, we present a novel algorithm to obtain diffeomorphic image or surface registrations with large deformations via quasi-conformal maps. The basic idea is to minimize an energy functional involving a Beltrami coefficient term, which measures the distortion of the quasi-conformal map. The Beltrami coefficient effectively controls the bijectivity and smoothness of the registration, even with very large deformations. Using the proposed algorithm, landmark-based registration between images or surfaces can be effectively computed. The obtained registration is guaranteed to be diffeomorphic (1-1 and onto), even with a large deformation or large number of landmark constraints. The proposed algorithm can also be combined with matching intensity (such as image intensity or surface curvature) to improve the accuracy of the registration. Experiments have been carried out on both synthetic and real data. Results demonstrate the efficacy of the proposed algorithm to obtain diffeomorphic registration between images or surfaces.
\end{abstract}

\begin{keywords}
Registration, large deformation, diffeomorphic, quasi-conformal mapping, Beltrami coefficient
\end{keywords}

%\begin{AMS}
%15A15, 15A09, 15A23
%\end{AMS}

\pagestyle{myheadings}
\thispagestyle{plain}
\markboth{Ka Chun Lam and Lok Ming Lui}{Diffeomorphic Registration with Large Deformations}

\section{Introduction}
Registration is a process of finding the optimal one-to-one correspondence between different data, such as images or surfaces. Applications can be found in various fields, including computer graphics, computer visions and medical imaging. For example, in medical imaging, finding accurate 1-1 correspondence between medical data is crucial for statistical shape analysis of the anatomical structures. While in computer graphics, surface registration is needed for texture mapping.

Different registration approaches have been developed. Existing algorithms can mainly be divided into three categories, namely, 1. landmark based registration, 2. intensity based registration and 3. hybrid registration using both landmark and intensity information. Landmark based registration computes a smooth 1-1 correspondence between corresponding data that matches important features. This kind of registration, with good feature alignment, is particularly crucial in medical imaging and computer graphics. For example, in computer graphics, landmark based registration is  used to obtain the constrained texture mapping. The main advantage of the landmark based method is that larger deformations can be dealt with and intuitive user-interaction can be incorporated. Intensity based registration aims to match corresponding data without feature landmarks. Registration is usually obtained by matching intensity functions, such as image intensity for image registration or surface curvature for surface geometric registration. The main advantage of the intensity based registration is that more image information is taken into account and the delineation of feature landmarks is not required. However, it usually cannot cope with large geometric deformations. Recently, hybrid registration that combines landmark based and intensity based methods have gained increased attention. Hybrid approaches use both the landmark and intensity information to guide the registration. This type of approaches can usually obtain more accurate registration result, since the advantages of landmark based and intensity based registration can be combined. In this work, we will mainly focus on the landmark based registration and the hybrid registration.

Most existing algorithms can compute registration accurately and efficiently when the deformation is small. However, the registration problem becomes challenging when large deformations occur. Bijectivity can be easily lost and overlaps can usually be observed in the obtained registration. This causes inaccuracies in the registration. It is therefore necessary to develop an algorithm to obtain diffeomorphic registration with large deformations.

In this paper, we introduce a novel method to obtain diffeomorphic image or surface registrations via quasi-conformal maps, which can deal with large deformations. The key idea is to minimize an energy functional involving a Beltrami coefficient term, which measures the distortion of the quasi-conformal map. The Beltrami coefficient effectively controls the bijectivity and smoothness of the registration, even with very large deformations. By minimizing the energy functional, we obtain an optimal Beltrami coefficient associated to the desired registration, which is guaranteed to be bijective. Using the proposed algorithm, landmark-based registration between images or surfaces can be effectively computed. The obtained registration is guaranteed to be diffeomorphic (1-1 and onto), even with a large deformation or a large number of landmark constraints. The proposed algorithm can also be combined with matching intensity (such as image intensity or surface curvature) to improve the accuracy of the registration. Numerical results show that the combination of landmark constraints with intensity matching can significantly improve the accuracy of the registration. To test the effectiveness of the proposed algorithm, experiments have been carried out on both synthetic and real data. Results show that the proposed algorithm can compute diffeomorphic registration between images or surfaces effectively and efficiently.

In summary, the contributions of this paper are three-folded. Firstly, we propose a variational method to search for an optimized Beltrami coefficient associated to a diffeomorphic quasi-conformal map with large deformations, which minimizes the local geometric distortion. Secondly, we apply the model to compute the landmark based registration, which can deal with very large deformations and large amount of landmark constraints. Thirdly, we extend the landmark based registration model to a hybrid registration model, which combine both landmark and intensity information to obtain more accurate registration.

The rest of the paper is organized as follows. In section \ref{previous}, we review some previous works closely related to this paper. In section \ref{math}, we describe some basic mathematical background related to our proposed model. In section \ref{model}, our proposed model for diffeomorphic registration with large deformations is explained in details. We describe the numerical implementation of the proposed algorithm in section \ref{numerical}. Experimental results are reported in section \ref{experiment}. Finally, we conclude our paper in section \ref{conclusion}.

\section{Previous works}\label{previous}
In this section, we will review some related works closely related to this paper.

Intensity-based image registration has been widely studied. A comprehensive survey on the existing intensity-based image registration can be found in \cite{SurveyZitova}. One of the commonly used method is based on the variational approaches to minimize the intensity mismatching error. For example, Vercauteren et al. \cite{Vercauteren} proposed the diffeomorphic demons registration algorithm, which is a non-parametric diffeomorphic image registration algorithm based on Thirion's demons algorithm\cite{Thirion}. The basic idea is to adapt the optimization procedure underlying the demons algorithm to a space of diffeomorphic transformations. The obtained registration is smooth and bijective. Several algorithms for surface registration that matches geometric quantities, such as curvatures, have also been propsoed \cite{Fischl2}\cite{Lyttelton}\cite{Lord}\cite{Yeo}. For example, Lyttelton et al. \cite{Lyttelton} proposed an algorithm for surface parameterizations based on matching surface curvatures. Yeo et al. \cite{Yeo} proposed the spherical demons method, which adopted the diffeomorphic demons algorithm \cite{Vercauteren}, to drive surfaces into correspondence based on the mean curvature and average convexity. Conformal surface registration, which minimizes angular distortions, has also been widely used to obtain a smooth 1-1 correspondence between surfaces \cite{Haker, Gu1, Gu3, Gu2, Hurdal, JinRicci, YLYangRicci, LuiVSRicci}. An advantage of conformal registrations is that they preserve local geometry well. Quasi-conformal surface registrations, which allows bounded amount of conformality distortion, have also been studied \cite{LuiBHF,LuiBHFHP,LuiQuasiYamabe,LuiBeltramirepresentation}. For example, Lui et al. \cite{LuiBHFHP} proposed to compute quasi-conformal registration between hippocampal surfaces based on the holomorphic Beltrami flow method, which matches geometric quantities (such as curvatures) and minimizes the conformality distortion \cite{LuiBHF}.

Landmark-based registration has also been widely studied and different algorithms have been proposed. Bookstein et al. \cite{Bookstein} proposed to use a thin-plate spline regularization (or biharmonic regularization) to obtain a registration that matches landmarks as much as possible. Tosun et al. \cite{Tosun} proposed to combine iterative closest point registration, parametric relaxation and inverse stereographic projection to align cortical sulci across brain surfaces. These diffeomorphisms obtained can better match landmark features, although not perfectly. Wang et al. \cite{Wang05,Luilandmark,Lui,Lui10} proposed to compute the optimized harmonic registrations of brain cortical surfaces by minimizing a compounded energy involving the landmark-mistmatching term \cite{Wang05,Luilandmark}. The obtained registration obtains an optimized harmonic map that better aligns the landmarks. However, landmarks cannot be perfectly matched, and bijectivity cannot be guaranteed under large number of landmark constraints. Later, Lin et al. \cite{Lin} propose a unified variational approach for registration of gene expression data to neuroanatomical mouse atlas in two dimensions that matches feature landmarks. Again, landmarks cannot be exactly matched. Inexact landmark-matching registrations are sometimes advantegous. In the case when landmark points/curves cannot be accurately delineated, this method is more tolerant of errors in labeling landmarks and gives better parameterization. In the situation when exact landmark matching is required, smooth vector field has been applied to obtain surface registration. Lui et al. \cite{Lui,Lui10} proposed the use of vector fields to represent surface maps and reconstruct them through integral flow equations. They obtained shape-based landmark matching harmonic maps by looking for the best vector fields minimizing a shape energy. The use of vector fields to compute the registration makes the optimization easier, although it cannot describe all surface maps. An advantage of this method is that exact landmark matching can be guaranteed. Time dependent vector fields can also be used \cite{Joshi,Glaunes,Glaunes2,Glaunes3,Glaunes4}. For example,  Glaun\'es et al. \cite{Glaunes} proposed to generate large deformation diffeomorphisms of a sphere, with given displacements of a finite set of template landmarks. The time dependent vector fields facilitate the optimization procedure, although it may not be a good representation of surface maps since it requires more memory. The computational cost of the algorithm is also expensive. Quasi-conformal mapping that matches landmarks consistently has also been proposed. Wei et al. \cite{Weiface} also proposed to compute quasi-conformal mappings for feature matching face registration. The Beltrami coefficient associated to a landmark points matching parameterization is approximated. However, either exact landmark matching or the bijectivity of the mapping cannot be guaranteed, especially when very large deformations occur.

Algorithms for hybrid registration, which combines both the landmark and intensity information to guide the registration, has also been proposed\cite{Hybrid1}\cite{Hybrid2}\cite{HybridHuang}\cite{Johnson02consistentlandmark}. For example, Christensen et al. \cite{Johnson02consistentlandmark} propsoed an algorithm for hybrid registration that uses both landmark and intensity information to guide the registration. The method utilizes the unidirectional landmark thin-plate spline (UL-TPS) registration technique together with a minimization scheme for the intensity difference to obtain good correspondence between images. Paquin et al. \cite{Hybrid1} proposed a registration method using a hybrid combination of coarse-scale landmark and B-splines deformable registration techniques. Chanwimaluang et al. \cite{Hybrid2} proposed a hybrid retinal image registration approach that combines both area-based and feature-based methods. Existing hybrid registration techniques can drive data into good correspondence when deformations are not too large. In this work, we propose a hybrid quasi-conformal registration method, which can deal with very large deformations.

\section{Mathematical background}\label{math}
In this work, we apply quasi-conformal maps to obtain diffeomorphic registrations with large deformations. In this section, we describe some basic theories related to quasi-conformal geometry. For details, we refer to readers to \cite{Gardiner}\cite{Lehto}.

A surface $S$ with a conformal structure is called a \emph{Riemann surface}. Given two Riemann surfaces $M$ and $N$, a map $f:M\to N$ is \emph{conformal} if it preserves the surface metric up to a multiplicative factor called the {\it conformal factor}. An immediate consequence is that every conformal map preserves angles. With the angle-preserving property, a conformal map effectively preserves the local geometry of the surface structure. %\cite{PDE:JLee}\cite{Add:Chern}\cite{Add:DOCARMO}.
\begin{figure}[t]
\centering
\includegraphics[height=1.35in]{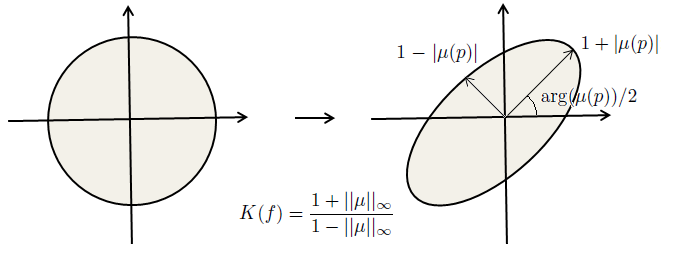}
\caption{Illustration of how the Beltrami coefficient determines the conformality distortion. \label{fig:illustration1}}
\end{figure}
A generalization of conformal maps is the \emph{quasi-conformal} maps, which are orientation preserving homeomorphisms between Riemann surfaces with bounded conformality distortion, in the sense that their first order approximations takes small circles to small ellipses of bounded eccentricity \cite{Gardiner}. Mathematically, $f \colon \mathbb{C} \to \mathbb{C}$ is quasi-conformal provided that it satisfies the Beltrami equation:
\begin{equation}\label{beltramieqt}
\frac{\partial f}{\partial \overline{z}} = \mu(z) \frac{\partial f}{\partial z}.
\end{equation}
\noindent for some complex-valued function $\mu$ satisfying $||\mu||_{\infty}< 1$. $\mu$ is called the \emph{Beltrami coefficient}, which is a measure of non-conformality. It measures how far the map at each point is deviated from a conformal map. In particular, the map $f$ is conformal around a small neighborhood of $p$ when $\mu(p) = 0$. Infinitesimally, around a point $p$, $f$ may be expressed with respect to its local parameter as follows:
\begin{equation}
\begin{split}
f(z) & = f(p) + f_{z}(p)z + f_{\overline{z}}(p)\overline{z} \\
& = f(p) + f_{z}(p)(z + \mu(p)\overline{z}).
\end{split}
\end{equation}

Obviously, $f$ is not conformal if and only if $\mu(p)\neq 0$. Inside the local parameter domain, $f$ may be considered as a map composed of a translation to $f(p)$ together with a stretch map $S(z)=z + \mu(p)\overline{z}$, which is postcomposed by a multiplication of $f_z(p),$ which is conformal. All the conformal distortion of $S(z)$ is caused by $\mu(p)$. $S(z)$ is the map that causes $f$ to map a small circle to a small ellipse. From $\mu(p)$, we can determine the angles of the directions of maximal magnification and shrinking and the amount of them as well. Specifically, the angle of maximal magnification is $\arg(\mu(p))/2$ with magnifying factor $1+|\mu(p)|$; The angle of maximal shrinking is the orthogonal angle $(\arg(\mu(p)) -\pi)/2$ with shrinking factor $1-|\mu(p)|$. Thus, the Beltrami coefficient $\mu$ gives us all the information about the properties of the map (See Figure \ref{fig:illustration1}).

The maximal dilation of $f$ is given by:
\begin{equation}
K(f) = \frac{1+||\mu||_{\infty}}{1-||\mu||_{\infty}}.
\end{equation}

Quasiconformal mapping between two Riemann surfaces $S_1$ and $S_2$ can also be defined. Instead of the Beltrami coefficient, the {\it Beltrami differential} is used. A Beltrami differential $\mu(z) \frac{\overline{dz}}{dz}$ on a Riemann surface $S$ is an assignment to each chart $(U_{\alpha},\phi_{\alpha})$ of an $L_{\infty}$ complex-valued function $\mu_{\alpha}$, defined on local parameter $z_{\alpha}$ such that
\begin{equation}
\mu_{\alpha}(z_{\alpha})\frac{d\overline{z_{\alpha}}}{dz_{\alpha}} = \mu_{\beta}(z_{\beta})\frac{d\overline{z_{\beta}}}{dz_{\beta}},
\end{equation}
\noindent on the domain which is also covered by another chart $(U_{\beta},\phi_{\beta})$. Here, $\frac{dz_{\beta}}{dz_{\alpha}}= \frac{d}{dz_{\alpha}}\phi_{\alpha \beta}$ and $\phi_{\alpha \beta} = \phi_{\beta}\circ \phi_{\alpha}^{-1}$.

\begin{figure*}[t]
\centering
\includegraphics[height=1.35in]{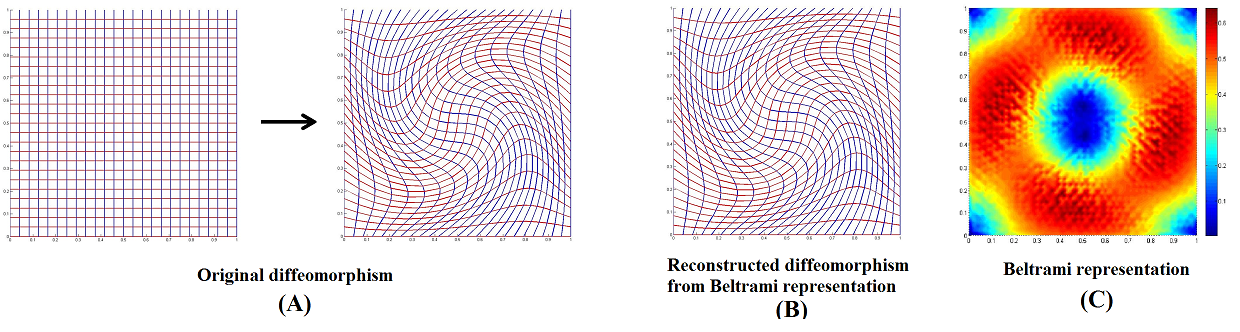}
\caption{Reconstruction of a diffeomorphism from its associated Beltrami coefficient. (A) shows a diffeomorphism between two rectangles. Its Beltrami coefficient is computed. (B) shows the reconstructed quasi-conformal map from the Beltrami coefficient. (C) shows the norm of the Beltrami coefficient.\label{fig:bcreconstruction}}
\end{figure*}
An orientation preserving diffeomorphism $f:S_1 \to S_2$ is called quasi-conformal associated with $\mu(z) \frac{\overline{dz}}{dz}$ if for any chart $(U_{\alpha},\phi_{\alpha})$ on $S_1$ and any chart $(V_{\beta},\psi_{\beta})$ on $S_2$, the mapping $f_{\alpha \beta}:= \psi_{\beta}\circ f\circ f_{\alpha}^{-1}$ is quasi-conformal associated with $\mu_{\alpha}(z_{\alpha})\frac{d\overline{z_{\alpha}}}{dz_{\alpha}}$.

\section{Proposed algorithm}\label{model}
In this section, we explain our proposed model for diffeomorphic registration with large deformation in details. The basic idea is to look for a quasi-conformal map to register two corresponding data, which can either be images or surfaces. The quasi-conformal map is obtained by minimizing an energy functional involving a Beltrami coefficient term, which measures the distortion of the quasi-conformal map. The Beltrami coefficient effectively controls the bijectivity and smoothness of the registration, even with very large deformations.

\subsection{Proposed model}

Let $S_1$ and $S_2$ be two corresponding images or surfaces. Our goal is to find a smooth 1-1 mapping $f:S_1\to S_2$ between $S_1$ and $S_2$ satisfying certain prescribed criteria. For landmark-based registration, we look for a registration that matches corresponding feature landmarks. Let $\{p_i \in S_1 \}_{i=1}^m$ and $\{q_i \in S_1 \}_{i=1}^m$ be the sets of corresponding feature landmarks defined on $S_1$ and $S_2$ respectively. We search for a diffeomorphism $f:S_1\to S_2$ subject to the landmark constraints that $f(p_i) = q_i$ for all $1 \leq i\leq m$.

We propose a variational approach to obtain an optimized quasi-conformal map $f$, which minimizes an energy functional $E_{LM}$ involving the Beltrami coefficient terms. More specifically, we propose to solve the following minimization problem:
\begin{equation}\label{landmarkmodel}
\begin{split}
f & = \mathbf{argmin}_{g:S_1\to S_2} E_{LM}(g) \\
&:= \mathbf{argmin}_{g:S_1\to S_2}\{\int_{S_1} |\nabla \mu_g|^2 + \alpha \int_{S_1} |\mu_g|^p\}
\end{split}
\end{equation}
\noindent subject to the constraints that:

\noindent C(i) $f(p_i) = q_i$ for $1\leq i\leq m$ (landmark constraint);

\noindent C(ii) $||\mu_f||_{\infty} <1$ (bijectivity).

The first term of $E_{LM}$ ensures the smoothness of $f$. The second term of $E_{LM}$ aims to minimize the conformality distortion of $f$. The constraint C(i) is the landmark constraint, which enforces $f$ to match corresponding landmarks consistently.

\begin{proposition}
If $f:S_1\to S_2$ satisfies the constraint C(ii), then $f$ is a diffeomorphism.
\end{proposition}
\begin{proof}
Suppose $f = u+ iv$ under some local coordinates. The Beltrami coefficient $\mu_f$ is given by:
\begin{equation}
\mu_f = \frac{\partial f}{\partial \overline{z}}/\frac{\partial f}{\partial z}
\end{equation}
where
\begin{equation}
\begin{split}
\frac{\partial f}{\partial \overline{z}} = (u_x - v_y) + i (u_y + v_x);\ \ \frac{\partial f}{\partial z} = (u_x + v_y) + i (v_x + u_y);
\end{split}
\end{equation}

Now, the Jacobian of $f$, $J_f$, is given by:
\begin{equation}
\begin{split}
& J_f = u_x v_y - u_y v_x\\
& = \frac{(u_x + v_y)^2 + (v_x + u_y)^2 - (u_x - v_y)^2 - (u_y + v_x)^2}{4}\\
& = |\frac{\partial f}{\partial z}|^2 -  |\frac{\partial f}{\partial \overline{z}} |^2 = |\frac{\partial f}{\partial z}|^2  (1 - |\mu_f|^2)
\end{split}
\end{equation}

Since $||\mu_f||_{\infty} < 1$, $|\frac{\partial f}{\partial z}|^2 \neq 0$. Also, $(1 - |\mu_f|^2) > 0$. Hence, $J_f > 0$ everywhere.

Since the Jacobian is postive everywhere, by the inverse function theorem, the mapping $f$ is locally invertible everywhere. In other words, $f$ is a diffeomorphism.
\end{proof}

\begin{proposition}[Landmark-matching registration]
Let:
\begin{equation}
\begin{split}
\mathcal{A} = & \{\nu \in C^1(\Omega_1): ||D\nu||_{\infty} \leq C_1; |\nu(p_i)|\leq C_2 + \epsilon; f^{\nu} (p_i) = q_i \mathrm{\ for\ } 1\leq i\leq n \}.
\end{split}
\end{equation}
Then: $E_{LM}$ has a minimizer in $\mathcal{A}$. In fact, $\mathcal{A}$ is compact.
\end{proposition}

\medskip

\begin{proof}
Note that $\mathcal{A} \neq \emptyset$. In particular, the unique Teichm\"uller map is indeed in $\mathcal{A}$. We first prove that $\mathcal{A}$ is complete. Let $\{\nu_n\}_{n=1}^{\infty}$ be a Cauchy sequence in $\mathcal{A}$ under the norm $||\cdot||_s$. Then, $\{\nu_n'\}_{n=1}^{\infty}$ is also a Cauchy sequence with respect to the $||\cdot||_{\infty}$ norm. Since $\mathcal{A}\subset H^1{\Omega_1}$, which is complete.  Thus, $\nu_n' \to g$ uniformly for some $g \in H^1(\Omega_1)$. Since $\nu_n'$ is continuous, $g$ is also continuous. Besides, $\nu_n (p_i)$  is convergent. Hence, $\nu_n\to \nu$ and $\nu_n \to \nu'$ uniformly for some $\nu \in C^1(\Omega_1)$.

In addition, $|\nu_n (p_i)| < C_1 + \epsilon$ for all $n$ implies $|\nu (p_i)| \leq C_1 + \epsilon$. $||D\nu_n||_{\infty}\leq C_2$ implies $||D\nu||_{\infty}\leq C_2$. Since $\nu_n \to \nu$ uniformly, $f^{v_n}\to f^{\nu}$ locally uniformly. This implies $f^{\nu}(p_i) = q_i$ for $1\leq i\leq n$. Therefore, $\nu \in \mathcal{A}$. Obviously, $\mathcal{A}$ is totally bounded. We conclude that $\mathcal{A}$ is compact.

Since $\mathcal{A}$ is compact and $E_{LM}$ is continuous in $\mathcal{A}$, $E_{LM}$ has a minimizer in $\mathcal{A}$.
\end{proof}

In order to improve the accuracy of the registration, one can combine the landmark-matching registration model with the intensity matching model. The intensities are functions defined on $S_1$ and $S_2$. Usually, they are image intensities for image registration and surface curvatures for surface registration. Ideally, we want to obtain a landmark-matching diffeomorphism $f:S_1\to S_2$ that matches the intensities as much as possible. We denote the intensities on $S_1$ and $S_2$ by $I_1:S_1\to \mathbb{R}$ and $I_2:S_2\to \mathbb{R}$ respectively. Our registration model can be modified as solving the following minimization problem:

\begin{equation}
\begin{split}
f & = \mathbf{argmin}_{g:S_1\to S_2} E_{IM}(g) \\
&:= \mathbf{argmin}_{g:S_1\to S_2}\{\int_{S_1} |\nabla \mu_g|^2 + \alpha \int_{S_1} |\mu_g|^p + \beta \int_{S_1} (I_1 - I_2(f))^2\}
\end{split}
\end{equation}
\noindent subject to the constraints C(i) and C(ii).

\begin{proposition}[Landmark and intensity matching registration]
$E_{IM}$ has a minimizer in $\mathcal{A}$ if $I_1$ and $I_2$ are continuous.
\end{proposition}
\begin{proof}
The exisitence of minimizer depends on the continuity of $I_1$ and $I_2$. If $I_1$ and $I_2$ are continuous, $E_{IM}$ is continuous. Since $\mathcal{A}$ is compact, $E_{IM}$ has a minimizer in $\mathcal{A}$.
\end{proof}

\subsection{Energy minimization}
In this subsection, we describe an algorithm to approximate the solutions of the above minimization problems.

\subsubsection{Landmark based registration model} Given two corresponding sets of landmarks $\{p_i\}_{i=1}^n$ and $\{q_i\}_{i=1}^n$on $S_1$ and $S_2$ respectively, our goal is to look for a diffeomorphism $f:S_1\to S_2$ that satisfies $f(p_i) = q_i$ ($i=1,...,n$) while minimizing the local geometric distortion. Our proposed model is to solve the variational problem (\ref{landmarkmodel}) as described in the last subsection.

More specifically, our goal is to look for an optimal Beltrami coefficient $\nu: S_1\to \mathbb{C}$, which is the Beltrami coefficient of some diffeomorphism $f:S_1\to S_2$, minimizing the following energy functional $E_{LM}$:
\begin{equation}
E_{LM}(\nu) = \int_{S_1} |\nabla \nu|^2 + \alpha \int_{S_1} |\nu|^p
\end{equation}
subject to the constraints that $||\nu||_{\infty}<1$, $f(p_i)=q_i$ for $i=1,2,...n$ and $\nu = \mu(f)$, where $\mu(f)$ is the Beltrami coefficient of $f$.

We apply a splitting method to solve the constrained optimization problem. In particular, we consider to minimize:
\begin{equation}
E_{LM}^{split} (\nu, f) = \int_{S_1} |\nabla \nu|^2 + \alpha \int_{S_1} |\nu|^p + \gamma_n \int_{S_1} |\nu -\mu(f)|^2
\end{equation}
subject to the constraints that $||\nu||_{\infty}<1$ and $f(p_i)=q_i$ for $i=1,2,...n$.

We iteratively minimize $E_{LM}^{split}$ subject to the constraints. Set $\nu_0 = 0$. Suppose $\nu_n$ is obtained at the $n^{th}$ iteration. Fixing $\nu_n$, we minimize $E_{LM}^{split}(\nu_n,f)$ over $f$, subject to the constraint that $f(p_i)=q_i$ ($i=1,2,...n$), to obtain $f_n$. Once $f_n$ is obtained, fixing $f_n$, we minimize $E_{LM}^{split}(\nu,f_n)$ over $\nu$ to obtain $\nu_{n+1}$.

To minimize $E_{LM}^{split}(\nu_n,f)$ over $f$ fixing $\nu_n$, it is equivalent to finding a landmark matching diffeomorphism $f_n:S_1\to S_2$, whose Beltrami coefficient closely resembles to $\nu_n$ and satisfies the landmark constraints $f(p_i) = q_i)$. To obtain such $f_n$, we propose to use the {\it Linear Beltrami Solver (LBS)} to find the descent direction for the Beltrami coefficients $\mu(f)$ such that it approaches to $\nu_n$ and the corresponding $f$ satisfies the landmark constraints.

Let $f = u + i v$. From the Beltrami equation (\ref{beltramieqt}),
\begin{equation}
\mu(f) = \frac{(u_x - v_y) + i\ (v_x + u_y)}{(u_x + v_y) + i(v_x - u_y)}
\end{equation}

Let $\mu(f) = \rho + i\ \tau $. We can write $v_x$ and $v_y$ as linear combinations of $u_x$ and $u_y$,
\begin{equation}\label{eqt:linearB1cont}
\begin{split}
-v_y & = \alpha_1 u_x + \alpha_2 u_y;\\
v_x & = \alpha_2 u_x + \alpha_3 u_y.
\end{split}
\end{equation}
\noindent where $\alpha_1 = \frac{(\rho -1)^2 + \tau^2}{1-\rho^2 - \tau^2} $; $\alpha_2 = -\frac{2\tau}{1-\rho^2 - \tau^2} $; $\alpha_3 = \frac{1+2\rho+\rho^2 +\tau^2}{1-\rho^2 - \tau^2} $.

Similarly,
\begin{equation} \label{eqt:linearB2cont}
\begin{split}
-u_y & = \alpha_1 v_x + \alpha_2 v_y;\\
u_x & = \alpha_2 v_x + \alpha_3 v_y.
\end{split}
\end{equation}

Since $\nabla \cdot \left(\begin{array}{c}
-v_y\\
v_x \end{array}\right) = 0$, we obtain
\begin{equation}\label{eqt:BeltramiPDE}
\nabla \cdot \left(A \left(\begin{array}{c}
u_x\\
u_y \end{array}\right) \right) = 0\ \ \mathrm{and}\ \ \nabla \cdot \left(A \left(\begin{array}{c}
v_x\\
v_y \end{array}\right) \right) = 0
\end{equation}

\noindent where $A = \left( \begin{array}{cc}\alpha_1 & \alpha_2\\
\alpha_2 & \alpha_3 \end{array}\right)$.

In the discrete case, the elliptic PDEs (\ref{eqt:BeltramiPDE}) can be discretized into sparse positive definite linear systems. Given $\nu_n$ and the landmark constraints, one can solve the linear systems with the landmark constraints in the least square sense. A landmark matching quasi-conformal map $f_n$, whose Beltrami coefficient closely resembles to $\nu_n$, can then be obtained.

Once $f_n$ is obtained, we minimize $E_{LM}^{split}(\nu, f_n)$ over $\nu$ while fixing $f_n$. In other words, we look for $\nu_{n+1}$ minimizing:
\begin{equation}
\int_{S_1} |\nabla \nu|^2 + \alpha \int_{S_1} |\nu|^p + \gamma_n \int_{S_1} |\nu -\mu(f_n)|^2
\end{equation}

By considering the Euler-Lagrange equation, it is equivalent to solving:
\begin{equation}\label{ELMsplit2continuous}
(\Delta + 2\alpha I + 2\gamma_n I) \nu_{n+1} = \mu(f_n)
\end{equation}

In discrete case, equation (\ref{ELMsplit2continuous}) can be discretized into a sparse linear system and can be solved efficiently. However, by simply minimizing equation (\ref{ELMsplit2continuous}) does not ensure that the corresponding diffeomorphism $f_n$ of the resultant Beltrami coefficient $\nu_n$ satisfies the landmark constraints. Therefore, we use the Linear Beltrami Solver with input $\nu_n$ together with the landmark constraints to obtain a mapping with the corresponding Beltrami coefficients $\nu_{n_i}$ resemble to $\nu_n$ and satisfies the exact landmark matching requirement. Using $d = \nu_{n_i} - \nu_n$ as a descent direction, we update $\nu_n$ from the solution of the Equation (\ref{ELMsplit2continuous}) by $\nu_n \leftarrow \nu_n + td$ for some small $t$. This guarantees the resultant $\nu_n$ is smooth and the landmark mismatch decreases.

We keep the iteration going to obtain a sequence of pair $\{(\nu_n, f_n)\}_{i=1}^{\infty}$. The iteration stops when $|\nu_{n+1}- \nu_n|< \epsilon$ for some small threshold $\epsilon$. Theoretically, the conventional penalty method requires that $\gamma_n$ increases in each iterations. In practice, we set $\gamma_n$ to be a large enough constant and the algorithm gives satisfactory results.

In summary, the proposed landmark based registration model can be described as follows:

\medskip

\noindent $\mathbf{Algorithm\ 1:}$ {\it(Landmark based registration)}\\
\noindent $\mathbf{Input:}$ {\it Images or surfaces: $S_1$ and $S_2$; cooresponding landmark sets $\{p_i \in S_1 \}_{i=1}^m$ and $\{q_i \in S_1 \}_{i=1}^m$.}\\
\noindent $\mathbf{Output:}$ {\it Optimal Beltrami coefficient $\nu^*$ and the landmark matching registration $f^*:S_1\to S_2$}\\
\vspace{-3mm}
\begin{enumerate}
\item {\it  Set $\nu_0 = 0$. Use LBS to reconstruct $f_0$ from $\nu_0$ satisfying the landmark constraints.;}
\item {\it Given $\nu_n$. Fixing $\nu_n$, obtain $f_n$ by LBS satisfying the landmark constraints. Fixing $f_n$, obtain $\nu_{n+1}$ by solving: $\nu_{n+1} = \textbf{argmin}_{\nu} \int | \nabla \nu |^2 + \alpha \int  |\nu|^p +  \gamma_n \int |\nu - \mu(f_n)|^2$}
\item Use LBS to obtain $\tilde{f_n}$ from $\nu_n$ satisfying the landmark constraints. Obtain $d = \mu(\tilde{f_n}) - \nu_n$ and update $\nu_n \leftarrow \nu_n + td$ for some small $t$.
\item {\it If $||\nu_{n+1} - \nu_n|| \geq \epsilon$, continue. Otherwise, stop the iteration.}
\end{enumerate}

\subsubsection{Hybrid registration model} The propsoed landmark based registration model can also be combined with matching intensity (such as image intensity for image registration or surface curvature for surface registration) to improve the accuracy of the registration result. More specifically, our goal is to look for an optimal Beltrami coefficient $\nu: S_1\to \mathbb{C}$, which is the Beltrami coefficient of some diffeomorphism $f:S_1\to S_2$, minimizing the following energy functional $E_{IM}$:
\begin{equation}
E_{IM}(\nu,f) = \int_{S_1} |\nabla \nu|^2 + \alpha \int_{S_1} |\nu|^p + \beta \int_{S_1} (I_1 - I_2(f))^2
\end{equation}
subject to the constraints that $||\nu||_{\infty}<1$ and $f(p_i) = q_i$ for $i=1,2,...,n$. Here, $I_1$ and $I_2$ are the intensity functions defined on $S_1$ and $S_2$ respectively.

We again apply a splitting method to solve the above constrained optimization problem. We consider to minimize:
\begin{equation}
\begin{split}
E_{IM}^{split} (\nu, \mu) &= \int_{S_1} |\nabla \nu|^2 + \alpha \int_{S_1} |\nu|^p + \sigma_n \int_{S_1} |\nu - \mu|^2\\
&+ \gamma \int_{S_1} (I_1 - I_2(f^{\mu}))^2
\end{split}
\end{equation}
subject to the constraints that $||\nu||_{\infty}<1$ and $f^{\mu}$ is the quasi-conformal map with Beltrami coefficient $\mu$ satisfying $f^{\mu}(p_i)=q_i$ for $i=1,2,...n$.

To solve the above optimization problem, we iteratively minimize $E_{IM}^{split}$ subject to the constraints. Set $\nu_0 = 0$. Suppose $\nu_n$ and $\mu_n$ is obtained at the $n^{th}$ iteration. Fixing $\nu_n$, we minimize $E_{IM}^{split}(\nu_n,\mu)$ over $\mu$, subject to the constraint that $f^{\mu}(p_i)=q_i$ ($i=1,2,...n$), to obtain $\mu_{n+1}$. Once $\mu_{n+1}$ is obtained, fixing $\mu_{n+1}$, we minimize $E_{IM}^{split}(\nu,\mu_{n+1})$ over $\nu$ to obtain $\nu_{n+1}$.

We first discuss the minimization $E_{IM}^{split}(\nu_n,\mu)$ over $\mu$ fixing $\nu_n$, subject to the constraint that $f^{\mu}(p_i)=q_i$ ($i=1,2,...n$). This problem is equivalent to solving:
\begin{equation}\label{IMsplit1}
\mu_{n+1} = \textbf{argmin}_{\mu} \gamma \int_{S_1} (I_1 - I_2(f^{\mu}))^2 + \sigma_n \int_{S_1}| \mu - \nu_n |^2
\end{equation}

Using the gradient descent method, we compute the descent direction $df$, which minimizes $\int_{S_1} (I_1 - I_2(f^{\mu}))^2$. $df$ is given by:
\begin{equation}
df = 2(I_1 - I_2(f^{\mu})) \nabla f^{\mu}.
\end{equation}

As $f^{\mu}$ is perturbed, its associated Beltrami coefficient is also adjusted by $d\mu_1$. The adjustment can be explicitly computed. Note that :
\begin{equation}
\frac{\partial (f + df)}{\partial \bar{z}} = (\mu + d \mu_1) \frac{\partial (f + df)}{\partial z}
\end{equation}
which implies:
\begin{equation}
\frac{\partial f}{\partial \bar{z}} + \frac{\partial df}{\partial \bar{z}} = \mu \frac{\partial f}{\partial z} + d \mu_1 \frac{\partial df}{\partial z} + \mu \frac{\partial df}{\partial z} + d \mu_1 \frac{\partial df}{\partial z}
\end{equation}

Thus, the adjustment $d\mu_1$ can be obtained by:
\begin{equation}
d \mu_1 = \left( \frac{\partial df}{\partial \bar{z}} - \mu \frac{\partial df}{\partial z}\right)/\frac{\partial f}{\partial z}
\end{equation}

Similiarly, we can obtain the descent direction $d\mu^2$, minimizes $\int_{S_1} |\mu-\nu_n|^2$. $d\mu_2$ is given by:
\begin{equation}
d\mu_2 = -2(\mu-\nu_n).
\end{equation}

Therefore, the descent direction to solve the optimization problem (\ref{IMsplit1}) is given by:
\begin{equation}
d\mu = \gamma d\mu_1 + \sigma_n d\mu_2
\end{equation}

Using the above formula for the descent direction, we obtain an updated Beltrami coefficient:
\begin{equation}
\tilde{\mu}_{n+1} = \mu_n + d\mu
\end{equation}

We then compute a quasi-conformal map $f_{n+1}$, whose Beltrami coefficient closely resembles to $\tilde{\mu}_{n+1}$, using LBS with the landmark constraints enforced. This step ensures a landmark matching registration can be obtained. We then update $\mu_n$ by: $\mu_{n+1} = \mu(f_{n+1})$.

Once $\mu_{n+1}$ is obtained, fixing $\mu_{n+1}$, we minimize $E_{IM}^{split}(\nu,\mu_{n+1})$ over $\nu$ to obtain $\nu_{n+1}$. In other words, we look for $\nu_{n+1}$ minimizing:
\begin{equation}
\int_{S_1} |\nabla \nu|^2 + \alpha \int_{S_1} |\nu|^p + \sigma_n \int_{S_1} |\nu -\mu_{n+1}|^2
\end{equation}

In the case when $p=2$, by considering the Euler-Lagrange equation, it is equivalent to solving:
\begin{equation}\label{ELMsplit2}
(\Delta + 2\alpha I + 2\sigma_n I) \nu_{n+1} = \mu_{n+1}
\end{equation}

In discrete case, equation (\ref{ELMsplit2}) can be discretized into a sparse linear system and can be solved efficiently. Similar to section 4.2.1, we use the Linear Beltrami Solver with input $\nu_n$ together with the landmark constraints to obtain a descent direction $d$ to update $\nu_n$ by $\nu_n \leftarrow \nu_n + td$ for some small $t$. This guarantees the resultant $\nu_n$ is smooth and the landmark mismatch decreases.

We keep the iteration going to obtain a sequence of pair $\{(\nu_n, \mu_n)\}_{i=1}^{\infty}$. The iteration stops when $|\mu_{n+1}- \mu_n|< \epsilon$ for some small threshold $\epsilon$. Again, the conventional penalty method requires that $\gamma_n$ increases in each iterations. However, in practice, we set $\sigma_n$ to be a large enough constant and the algorithm gives satisfactory results.

In summary, the proposed hybrid registration model can be described as follows:

\noindent $\mathbf{Algorithm\ 2:}$ {\it(Hybrid registration)}\\
\noindent $\mathbf{Input:}$ {\it Images or surfaces: $S_1$ and $S_2$; cooresponding landmark sets $\{p_i \in S_1 \}_{i=1}^m$ and $\{q_i \in S_1 \}_{i=1}^m$; intensity functions $I_1$ and $I_2$ defined on $S_1$ and $S_2$ respectively.}\\
\noindent $\mathbf{Output:}$ {\it Optimal Beltrami coefficient $\mu^*$ and the landmark matching registration $f^*:S_1\to S_2$}\\
\vspace{-3mm}
\begin{enumerate}
\item {\it  Set $\nu_0 = 0$. Use LBS to reconstruct $f_0$ from $\tilde{\mu}_0 = 0$ satisfying the landmark constraints. Set $\mu_0 = \mu(f_0)$;}
\item {\it Given $\nu_n$ and $\mu_n$. Fixing $\nu_n$, obtain $\tilde{\mu}_{n+1}$ by solving: $$ \tilde{\mu}_{n+1} = \textbf{argmin}_{\mu} \Big\{ \gamma \int_{S_1} (I_1 - I_2(f^{\mu}))^2 + \sigma_n \int_{S_1}| \mu - \nu_n |^2 \Big\}$$.}
\item {\it Using LBS, compute $f_{n+1}$ whose Beltrami coefficient closely resembles to $\tilde{\mu}_{n+1}$ with landmark constraints enforced. Let $\mu_{n+1} = \mu(f_{n+1})$.}
\item {\it Fixing $\mu_{n+1}$, obtain $\nu_{n+1}$ by solving:\\
$\nu_{n+1} = \textbf{argmin}_{\nu} \int | \nabla \nu |^2 + \alpha \int  |\nu|^p +  \sigma_n \int |\nu - \mu_{n+1})|^2$}

\item Use LBS to obtain $\tilde{f}_{n+1}$ from $\nu_{n+1}$ satisfying the landmark constraints. Obtain $d = \mu(\tilde{f}_{n+1}) - \nu_{n+1}$ and update $\nu_{n+1} \leftarrow \nu_{n+1} + td$ for some small $t$.

\item {\it If $||\mu_{n+1} - \mu_n|| \geq \epsilon$, continue. Otherwise, stop the iteration.}
\end{enumerate}

\section{Numerical implementation}\label{numerical}
The proposed models for landmark based and hybrid registration rely on the Linear Beltrami Solver(LBS) and solving the Euler-Lagrange(E-L) equations. In this section, we will describe the numerical implementation of the LBS and also the discretization of equations.

Practically speaking, 2D domains or surfaces in $\mathbb{R}^3$ are usually represented discretely by triangular meshes. Suppose $K_1$ and $K_2$ are two surface meshes with the same topology representing $S_1$ and $S_2$. We define the set of vertices on $K_1$ and $K_2$ by $V^1 = \{v_i^1\}_{i=1}^n$ and $V^2 = \{v_i^2\}_{i=1}^n$ respectively. Similarly, we define the set of triangular faces on $K_1$ and $K_2$ by $F^1 = \{T_j^1\}_{j=1}^m$ and $F^2 = \{T_j^2\}_{j=1}^m$.

\subsection{Numerical details of LBS}

Suppose $f:K_1 \to K_2$ is an orientation preserving piecewise linear homeomorphism between $K_1$ and $K_2$. We can assume $K_1$ and $K_2$ are both embedded in $\mathbb{R}^2$. In case $K_1$ and $K_2$ are surface meshes in $\mathbb{R}^3$, we first parameterize them conformally by $\phi_1:K_1\to D_1 \subseteq \mathbb{R}^2$ and $\phi_2:K_2\to D_2\subseteq \mathbb{R}^2$. The composition of $f$ with the conformal parameterizations, $\tilde{f}:= \phi_2 \circ f\circ\phi_1^{-1}$, is then an orientation preserving piecewise linear homeomorphism between $D_1$ and $D_2$ embedded in $\mathbb{R}^2$.

To compute the quasi-conformal mapping, the key idea is to discretize Equation \ref{eqt:BeltramiPDE} with two linear systems.

Given a map $f=(u+iv): K_1 \to K_2$, we can easily compute its associated Beltrami coefficient $\mu_f$, which is a complex-valued function defined on each triangular faces of $K_1$. To compute $\mu_f$, we simply need to approximate the partial derivatives on every face $T$. We denote them by $D_x f(T) = D_x u + i D_x v$ and  $D_y f(T) = D_y u + i D_y v$ respectively. Note that $f$ is piecewise linear. The restriction of $f$ on each triangular face $T$ can be written as:
\begin{equation}
f|_T (x,y) = \left( \begin{array}{c}
a_T x + b_T y + r_T \\
c_T x + d_T y + s_T \end{array} \right)
\end{equation}

Hence, $D_x u(T) = a_T$, $D_y u(T) = b_T$,  $D_x v(T) = c_T$ and $D_y v(T) = d_T$. Now, the gradient:
\begin{equation}
\nabla _T f := (D_x f(T), D_y f(T))^t
\end{equation}
on each face $T$ can be computed by solving the linear system:
\begin{equation}\label{eqt:gradient}
\left( \begin{array}{c}
\vec{v}_1 - \vec{v}_0\\
\vec{v}_2 - \vec{v}_0\end{array} \right)\nabla_T \tilde{f}_i = \left( \begin{array}{c}
\frac{\tilde{f}_i(\vec{v}_1) - \tilde{f}_i(\vec{v}_0)}{|\vec{v}_1 - \vec{v}_0|}\\
\frac{\tilde{f}_i(\vec{v}_2) - \tilde{f}_i(\vec{v}_0)}{|\vec{v}_2 - \vec{v}_0|}\end{array} \right),
\end{equation}
\noindent where $[\vec{v_0},\vec{v_1}]$ and $[\vec{v_0},\vec{v_2}]$ are two edges on $T$. By solving equation (\ref{eqt:gradient}), $a_T$, $b_T$, $c_T$ and $d_T$ can be obtained. The Beltrami coefficient $\mu_f(T)$ of the triangular face $T$ can then be computed from the Beltrami equation (\ref{beltramieqt}) by:
\begin{equation}\label{eqt:BC}
\mu_f(T) = \frac{(a_T - d_T)+i(c_T + b_T)}{(a_T + d_T)+i(c_T - b_T)},
\end{equation}

Equation (\ref{eqt:linearB1cont}) and (\ref{eqt:linearB2cont}) are both satisfied on every triangular faces. Let $\mu_f(T) = \rho_T + i\ \tau_T$. The discrete versions of Equation (\ref{eqt:linearB1cont}) and (\ref{eqt:linearB2cont}) can be obtained.
\begin{equation}\label{eqt:BCsplit1}
\begin{split}
-d_T & = \alpha_1(T) a_T + \alpha_2(T) b_T\\
 c_T & = \alpha_2(T) a_T + \alpha_3(T) b_T
\end{split}
\end{equation}

and

\begin{equation}\label{eqt:BCsplit2}
\begin{split}
-b_T & = \alpha_1(T) c_T + \alpha_2(T) d_T\\
 a_T & = \alpha_2(T) c_T + \alpha_3(T) d_T
\end{split}
\end{equation}
\noindent where:
$\alpha_1(T) = \frac{(\rho_T -1)^2 + \tau_T^2}{1-\rho_T^2 - \tau_T^2} $; $\alpha_2(T) = -\frac{2\tau_T}{1-\rho_T^2 - \tau_T^2} $ and

\noindent $\alpha_3(T) = \frac{1+2\rho_T+\rho_T^2 +\tau_T^2}{1-\rho_T^2 - \tau_T^2} $.

In order to discretize Equation (\ref{eqt:BeltramiPDE}), we need to introduce the {\it discrete divergence}. The discrete divergence can be defined as follows. Let $T = [v_i,v_j, v_k]$ and $w_I = f(v_I)$ where $I=i,j$ or $k$. Suppose $v_I = g_I + i\ h_I$ and $w_I = s_I + i\ t_I$ ($I=i,j,k$). Using equation (\ref{eqt:gradient}), $a_T, b_T, c_T$ and $d_T$ can be written as follows:
\begin{equation}
\begin{split}
a_T &= A_i^T s_i + A_j^T s_j + A_k^T s_k;\\
b_T &= B_i^T s_i + B_j^T s_j + B_k^T s_k;\\
c_T &= A_i^T t_i + A_j^T t_j + A_k^T t_k;\\
d_T &= B_i^T t_i + B_j^T t_j + B_k^T t_k;
\end{split}
\end{equation}
\noindent where:
\begin{equation}
\begin{split}
&A_i^T = (h_j-h_k )/Area(T),\\
&A_j^T = (h_k-h_i )/Area(T),\\
&A_k^T = (h_i-h_j )/Area(T);\\
&B_i^T = (g_k-g_j )/Area(T),\\
&B_j^T = (g_i-g_k )/Area(T),\\
&B_k^T = (g_j-g_i )/Area(T);
\end{split}
\end{equation}

Suppose $\vec{V} = (V_1, V_2)$ is a discrete vector field defined on every triangular faces. For each vertex $v_i$, let $N_i$ be the collection of neighborhood faces attached to $v_i$. We define the discrete divergence $Div $ of $\vec{V}$ as follows:
\begin{equation}
Div  (\vec{V})(v_i) = \sum_{T\in N_i} A_i^T V_1(T) + B_i^T V_2(T)
\end{equation}

By careful checking, one can prove that
\begin{equation}
\sum_{T\in N_i} A_i^T b_T = \sum_{T\in N_i} B_i^T a_T;\ \sum_{T\in N_i} A_i^T d_T = \sum_{T\in N_i} B_i^T c_T.
\end{equation}

This gives,
\begin{equation}
Div\ \left(\begin{array}{c}
-D_y u\\
D_x u \end{array}\right) = 0 \ \ \mathrm{and}\ \ Div \left(\begin{array}{c}
-D_y v\\
D_x v\end{array}\right) = 0
\end{equation}

As a result, Equation (\ref{eqt:BeltramiPDE}) can be discretized:
\begin{equation}\label{eqt:BeltramiPDEdiscrete}
Div \left( A \left(\begin{array}{c}
D_x u\\
D_y u\end{array}\right) \right) = 0\ \ \mathrm{and}\ \ Div \left(A \left(\begin{array}{c}
D_x v\\
D_y v\end{array}\right) \right) = 0
\end{equation}
\noindent where $A = \left(\begin{array}{cc}
\alpha_1 & \alpha_2\\
\alpha_2 & \alpha_3 \end{array}\right) $. This is equivalent to:

\begin{equation}\label{eqt:linearB12}
\begin{split}
\sum_{T\in N_i} A_i^T [\alpha_1(T) a_T  &+ \alpha_2(T) b_T] + B_i^T[\alpha_2(T) a_T + \alpha_3(T) b_T] = 0\\
\sum_{T\in N_i} A_i^T [\alpha_1(T) c_T &+ \alpha_2(T) d_T] + B_i^T[\alpha_2(T) c_T + \alpha_3(T) d_T] = 0
\end{split}
\end{equation}

\noindent for all vertices $v_i \in D$. Note that $a_T$, $b_T$, $c_T$ and $d_T$ can be written as a linear combination of the x-coordinates and y-coordinate of the desired quasi-conformal map $f$. Hence, equation (\ref{eqt:linearB12}) gives us the linear systems to solve for the x-coordinate and y-coordinate function of $f$.

Besides, $f$ has to satisfy certain constraints on the boundary. One common situation is to give the Dirichlet condition on the whole boundary. That is, for any $v_b \in \partial K_1$
\begin{equation}
f(v_b) = w_b \in \partial K_2
\end{equation}

Note that the Dirichlet condition is not required to be enforced on the whole boundary. The proposed algorithm also allows free boundary condition. For example, in the case that $K_1$ and $K_2$ are rectangles, the desired quasi-conformal map should satisfy
\begin{equation}\label{eqt:boundary1}
\begin{split}
&f(0) = 0; f(1) = 1\ f(i) = i\ f(1+i) = 1+i;\\
&\mathbf{Re}(f) = 0 \mathrm{\ on\ arc\ }[0, i];\ \mathbf{Re}(f) = 1 \mathrm{\ on\ arc\ }[1, 1+i];\\
&\mathbf{Imag}(f) = 0 \mathrm{\ on\ arc\ }[0, 1];\ \mathbf{Imag}(f) = 1 \mathrm{\ on\ arc\ }[i, 1+i]
\end{split}
\end{equation}

Besides, in our case, interior landmark correspondences $\{p_i\}_{i=1}^n \leftrightarrow \{q_i\}_{i=1}^n$ are also enforced. Thus, we should add this constraint to the linear systems. Mathematically, it is described as $f(p_i)= q_i$ ($i=1,2,...,n$).

\subsection{Discretization of E-L equations}

In discrete case, as the Beltrami coefficient $\mu(T)$ is defined on each triangular face, we first approximate $\mu(v)$ by

\begin{equation}\label{eqt:smoothmu1}
\mu(v_i) = \frac{1}{N_i} \sum_{T \in N_i} \mu(T)
\end{equation}

\noindent where $N_i$ is the collection of neighborhood faces attached to $v_i$. The obtained $\mu(v_i)$ will be the average of the Beltrami coefficients $\mu(T)$ on 1-ring neighbourhood triangles.

Discretizing Equation (\ref{ELMsplit2continuous}) turns down to the problem of discretizing the laplacian operator $\Delta$. Let $T_1 = [v_i, v_j, v_k]$ and $T_2 = [v_i, v_j, v_l]$. The mesh laplacian is defined to be
\begin{equation}\label{laplaciandiecrete}
\Delta(f(v_i)) = \sum_{T \in N_i} \frac{ \cot \alpha_{ij} + \cot \beta_{ij}}{2} (f(v_j) - f(v_i))
\end{equation}

\noindent where $\alpha_{ij}$ and $\beta_{ij}$ are the two interior angles of $T_1$ and $T_2$ which are opposite to the edge $[v_i, v_j]$. To find $\alpha_{ij}$ and similar $\beta_{ij}$, we follows the idea of \cite{ETH_laplacian}. Let $l_{ij}$ be the length of the edge $[v_i, v_j]$. By law of cosines: $l_{ij}^2 = l_{jk}^2 + l_{ki}^2 - 2l_{jk} l_{ki} \cos \alpha_{ij},$
\noindent we have
\begin{equation}\label{cosine1}
\cos \alpha_{ij} = \frac{-l_{ij}^2 + l_{jk}^2 + l_{ki}^2}{2l_{jk}l_{ki}}.
\end{equation}
\noindent Similar, by the law of sines: $\text{Area}(T_1) = \frac{1}{2} l_{jk} l_{ki} \sin \alpha_{ij}$, we have
\begin{equation}\label{sine1}
\sin \alpha_{ij} = \frac{2Area(T_1)}{l_{jk} l_{ki}}.
\end{equation}
\noindent Therefore we have
\begin{equation}\label{cotangent1}
\cot \alpha_{ij} = \frac{-l_{ij}^2 + l_{jk}^2 + l_{ki}^2}{4 \text{Area}(T_1)}
\end{equation}
\noindent and the discrete laplace beltrami operator can be constructed.

As the solution of the equation (\ref{ELMsplit2continuous}),  $\nu_{n+1}$ , is defined on vertices, we have to approximate $\nu_{n+1}(T)$ on each face before proceeding to next step . The approximation is taken to be

\begin{equation}\label{muv2f}
\nu_{n+1}(T) = \frac{1}{3} \sum_{v_i \in T} \nu_{n+1}(v_i)
\end{equation}

\subsection{Intensity matching}

In section 4.2.2, we propose to solve
$$\textbf{argmin}_{\mu}  \Big\{ \gamma \int_{S_1} (I_1 - I_2(f^{\mu}))^2 + \sigma_n \int_{S_1}| \mu - \nu_n |^2 \Big\}. $$
\noindent by using gradient descent. However, the steps involves calculations of the gradient of $\mu$, which is a second order derivatives. It may cause computation error and instability. We therefore separately search for the descent direction of $\mu_1$ and $\mu_2$ for the two terms $\int_{S_1} (I_1 - I_2(f^{\mu}))^2$ and $\int_{S_1} | \mu - \nu_n |^2 $ respectively in the following way.

For the first term, we apply the Demon force proposed by Wang et al. \cite{DemonForce} to find the deformation:

\begin{equation}\label{demonforce}
\begin{split}
u =& \frac{(I_1 - I_2)\nabla (I_2)}{| \nabla (I_2) |^2 + \alpha^2 ( I_1 - I_2)^2} + \frac{(I_1 - I_2)\nabla (I_1)}{| \nabla (I_1)|^2 + \alpha^2 ( I_1 - I_2)^2}
\end{split}
\end{equation}

\noindent where $u$ is the deformation vector field. The corresponding Beltrami coefficient of the deformation is

\begin{equation}\label{demonBC}
\mu_d = \frac{\partial (Id + u)}{\partial \bar{z}}\left/ \frac{\partial (Id +u)}{\partial z} \right.
\end{equation}

By the composition rule of the Beltrami coefficient, we have
\begin{equation}
\mu( u(f) ) = \frac{\mu(f) + \frac{\overline{f_z}}{f_z}\mu_d}{1 + \frac{\overline{f_z}}{f_z} \overline{\mu(f)}\mu_d}
\end{equation}

Then the descent direction of $d\mu_1$, minimizes $\int_{S_1} (I_1 - I_2(f^{\mu}))^2$, is approximated by
\begin{equation}
d\mu_1 \approx \mu( u(f) ) - \mu ( f )
\end{equation}

For the second term, we can obtain the descent direction as $d\mu_2 = -2(\mu-\nu_n)$. Therefore, the descent direction to solve the optimization problem is given by:
\begin{equation}
d\mu = \gamma d\mu_1 + \sigma_n d\mu_2
\end{equation}

Using the Demon force as registration guarantee the smoothness of $\mu_d$ and also stabilize the calculation of descent direction.

\subsection{Multiresolution scheme}
To reduce the computation cost of registering high resolution images (high quality surface meshes), we adopt a multiresolution scheme for the registration procedure. In the scheme, we first coarsen the both $I_1$ and $I_2$ by $k$ layers, where $I_j^0 = I_j$ and $I_j^k$ is the coarest images of $I_j$, for $j = 1,2$. Registration process is then carried out to register between $I_1^k$ and $I_2^k$. Diffeomorphism $f_k$ can then be obtained. To proceed to finer scale, we adopt a linear interpolation on $f_k$ to obtain $f_{k-1}$, which serves as the initial map for regstration in finer layer. Resultant mapping of this scheme will be the diffeomorphism which matches the intensity and satisfies the landmark constraints for the original image resolution.

\section{Experimental results}\label{experiment}
We have test our proposed algorithms on synthetic data together with real medical data. In this section, experimental results are reported.

\begin{figure*}[t]
\centering
\includegraphics[height=1.75in]{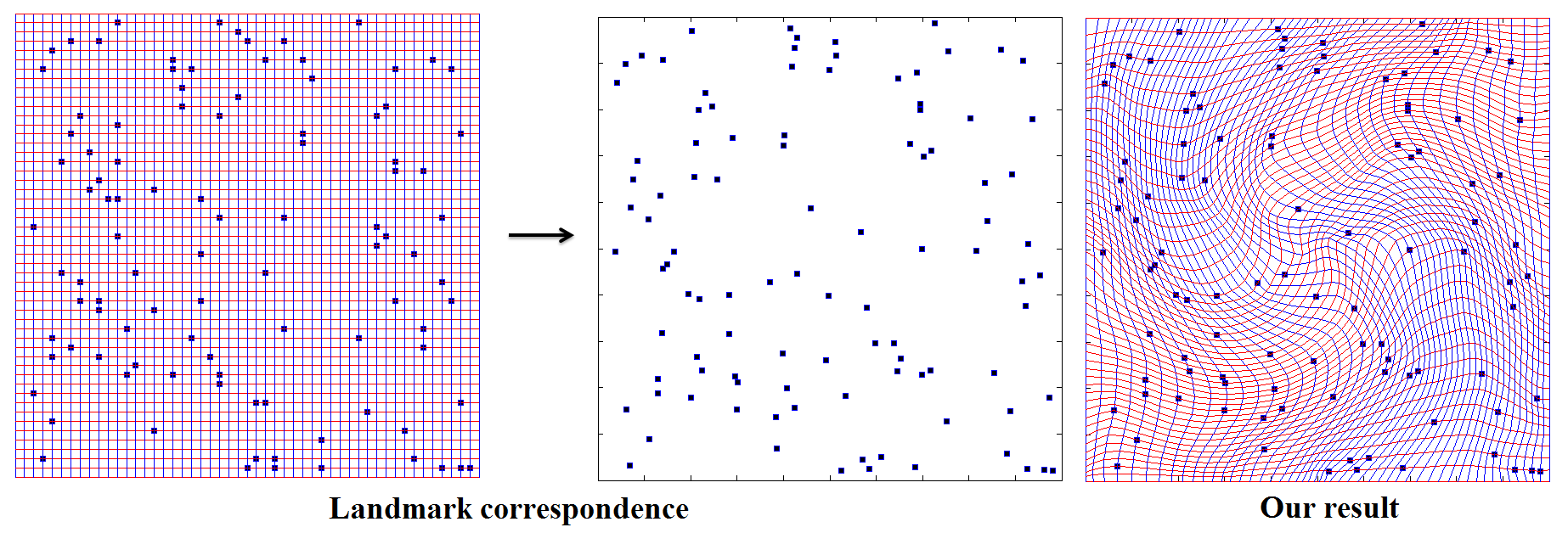}
\caption{Landmark based registration with large amount of landmark constraints. (A) shows the the correspondence between two landmark sets defined on two unit squares. (B) shows the obtained landmark matching diffeomorphic registration using our proposed algorithm. \label{fig:landmark_example2}}
\end{figure*}

\subsection{Landmark based registration}
\paragraph{Example 1} In this example, we test our proposed landmark based registration model on a synthetic example with large amount of landmark constraints enforced. Figure\ref{fig:landmark_example2}(A) shows the the correspondence between two landmark sets defined on two rectangles. 78 corresponding landmark features are labeled on each rectangles. We compute the landmark matching diffeomorphic registration between the two rectangles, using our proposed landmark based registration model. The registration result is as shown in (B), which is visualized by the deformation of the regular grids under the registration. Note that the obtained registration is bijective. No overlaps or flips can be found.

\begin{figure*}[t]
\centering
\includegraphics[height=1.75in]{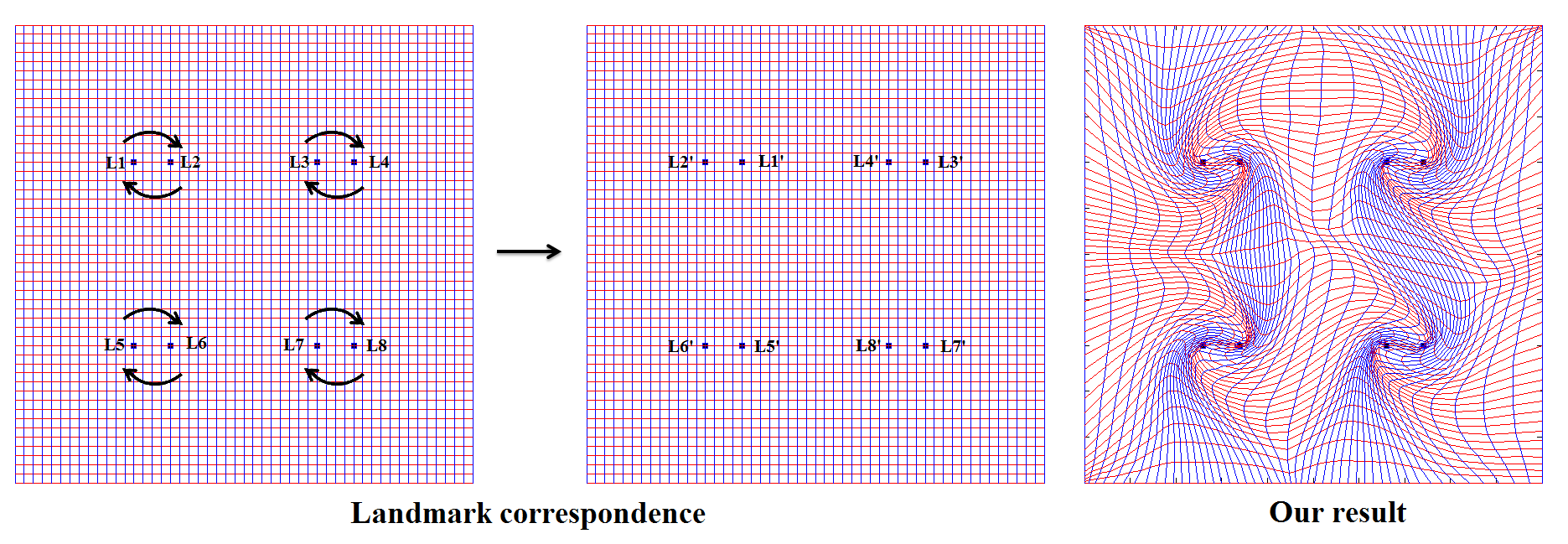}
\caption{Landmark based registration with large deformations. (A) shows the the correspondence between two landmark sets defined on two unit squares. (B) shows the obtained landmark matching diffeomorphic registration using our proposed algorithm. \label{fig:landmark_example1}}
\end{figure*}

\paragraph{Example 2} We first test our proposed algorithm on a synthetic example to obtain a landmark matching registration between two 2D rectangles with very large deformations. Figure \ref{fig:landmark_example1}(A) shows two rectangles, with corresponding landmark sets defined on each of them. The presecribed deformations of the landmarks are large. Using our proposed landmark based registration model, we compute the landmark matching diffeomorphic registration between the two rectangles. The registration result is as shown in (B), which is visualized by the deformation of the regular grids under the registration. Note that the obtained registration is bijective. No overlaps or flips can be found.

\begin{figure*}[t]
\centering
\includegraphics[height=3.35in]{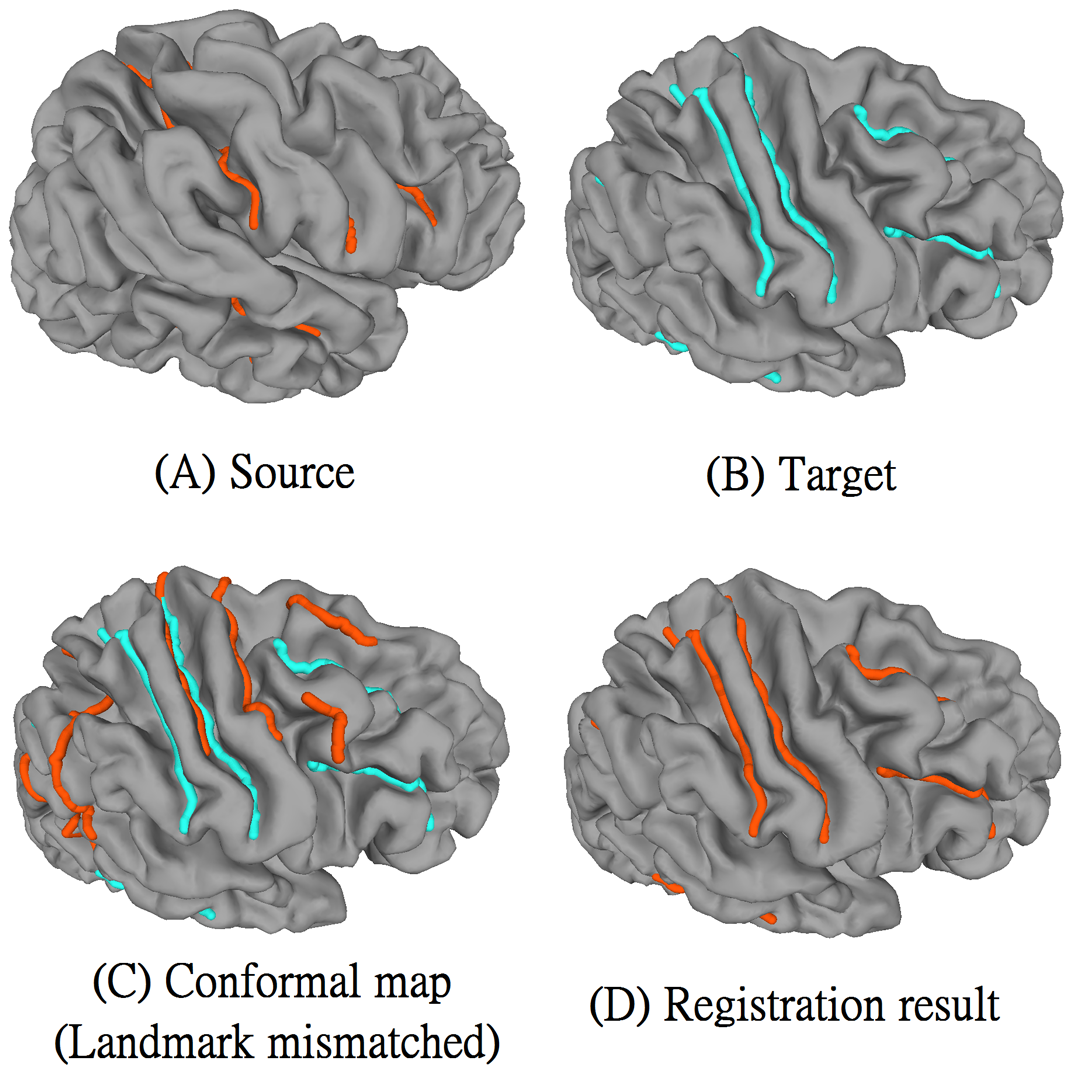}
\caption{(A) and (B) shows two brain cortical surfaces, each of them are labeled with six corresponding sulcal landmarks. (C) shows the conformal registration between the two surfaces without landmark matching. (D) shows the registration result using our proposed landmark-matching quasi-conformal registration.\label{fig:mistmatch}}
\end{figure*}

\begin{figure*}[t]
\centering
\includegraphics[height=1.75in]{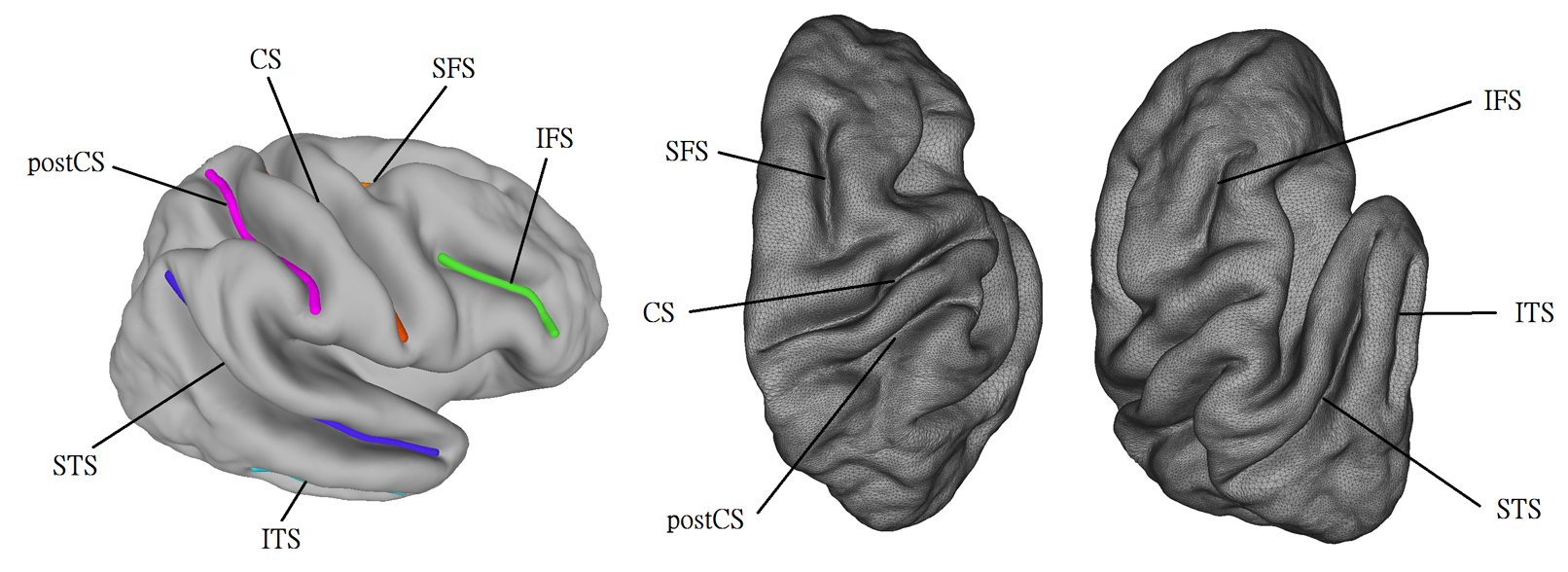}
\caption{Using the proposed algorithm, the landmark-matching quasi-conformal registrations with 6 sulcal landmarks between 10 brain cortical surfaces are computed. The mean surface is then computed after the cortical surfaces are registered. The sulcal features are well-preserved, illustrating that the landmarks are consistently matched under the registration.\label{fig:MeanSurface}}
\end{figure*}

\paragraph{Example 3} (Brain landmark matching registration) We apply the proposed algorithm to compute landmark matching quasi-conformal registration between brain cortical surfaces. Figure \ref{fig:mistmatch}(A) and (B) show two brain cortical surfaces, each of them are labeled by 6 sulcal landmarks. Using our proposed method, we compute the landmark-matching quasi-conformal registration between them. Figure \ref{fig:mistmatch}(C) shows the conformal registration between the two surfaces. Note that the corresponding landmarks cannot be matched. Figure \ref{fig:mistmatch}(D) shows the registration result using our proposed landmark-matching quasi-conformal registration, which matches corresponding landmarks consistently. Figure shows the energy plot versus iterations. It demonstrates our method iteratively minimizes the energy functional to an optimal quasi-conformal map between the two brain surfaces. In Figure \ref{fig:MeanSurface}, we compute the landmark-matching quasi-conformal registrations with 6 sulcal landmarks between 10 brain cortical surfaces. The mean surface is then computed after the cortical surfaces are registered. The sulcal features are well-preserved, illustrating that the landmarks are consistently matched under the registration.

\begin{figure*}[t]
\centering
\includegraphics[height=2.25in]{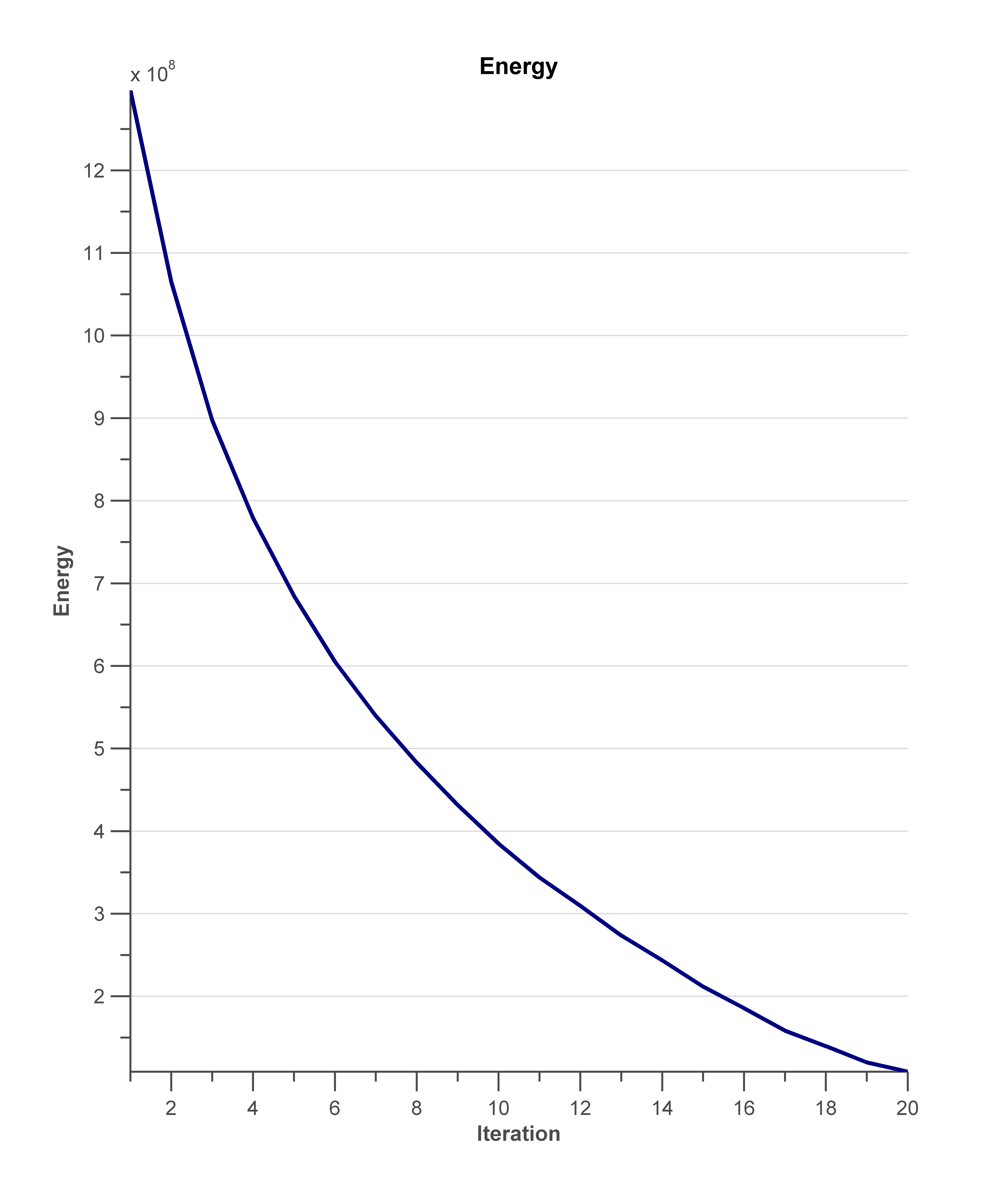}
\caption{The energy plot versus iterations for the landmark-matching quasi-conformal registration between brain cortical surfaces.\label{fig:Brainsurfaceenergyplot}}
\end{figure*}

\paragraph{Example 4} We also compared our proposed landmark-matching quasi-conformal registration algorithm with other state-of-the-art algorithms, namely, 1. harmonic map; 2. thin-plate spline (TPS) and 3. LDDMM. We compared our method with others with different sizes of deformation. As shown in Table \ref{table:comparisondeformation}, our method outperforms other methods. In all cases (tiny, moderate and large deformation), our method computes a non-overlap landmark-matching registration with the least amount of time. Both harmonic map and TPS has overlaps for their obtained landmark-matching registration results, although the computations of these methods are efficient. LDDMM can obtain non-overlapping landmark-matching registration results, however, the computational cost is comparatively much more expensive. 

\begin{table}[!h]
\caption{Comparison with other methods with different sizes of deformation}
\begin{center}
\begin{tabular}{c||c|c|c}
(Time / Overlap)& Tiny & Moderate & Large \\
\hline
QC & 6.220 s / 0 & 9.632 s / 0 & 12.934 s / 0 \\
Harmonic Map & 1.633 s / 13 & 1.665 s / 42 & 1.652 s / 110 \\
TPS & 0.308 s / 20 & 0.339 s / 27 & 0.253 s / 27 \\
LDDMM & 382.316 s / 0 & 396.240 s / 0 & 409.902 s / 0 \\
\end{tabular}
\label{table:comparisondeformation}
\end{center}
\end{table}

\subsection{Hybrid registration}

\begin{figure*}[t]
\centering
\includegraphics[height=1.75in]{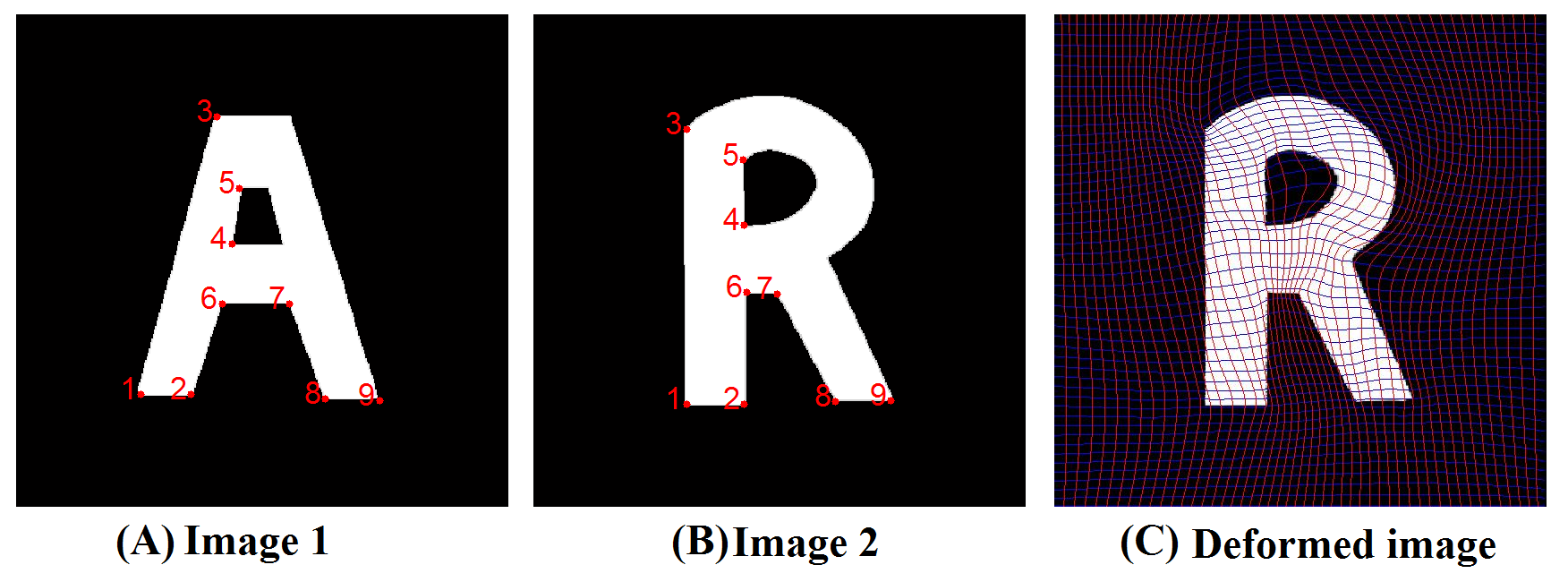}
\caption{(A) and (B) shows two images to be registered. Corresponding feature landmarks are labeled on each images. (B) shows the obtained diffeomorphic registration using our proposed hybrid registration model.\label{fig:AR_result}}
\end{figure*}

\begin{figure*}[t]
\centering
\includegraphics[height=2.35in]{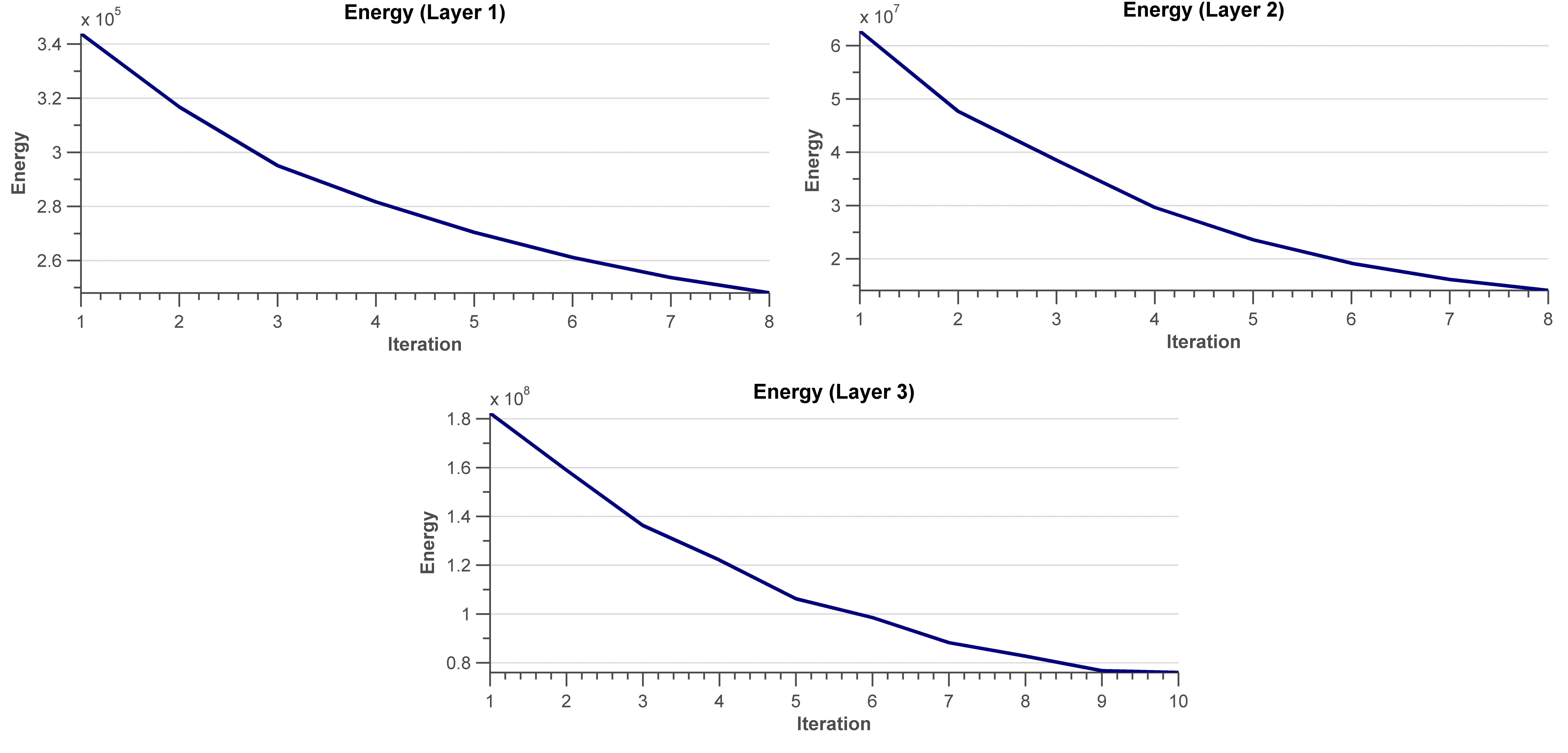}
\caption{The plots of energy versus iterations for the hybrid registration between the 'A' and 'R' images. Multi-level approach is applied to perform the registration from the coarsest layer to the finest (original resolution) layer. The energy plots at different layers are shown.\label{fig:ARenergy}}
\end{figure*}

\begin{figure*}[t]
\centering
\label{fig:multiresolutionAR}
\includegraphics[height=1.85in]{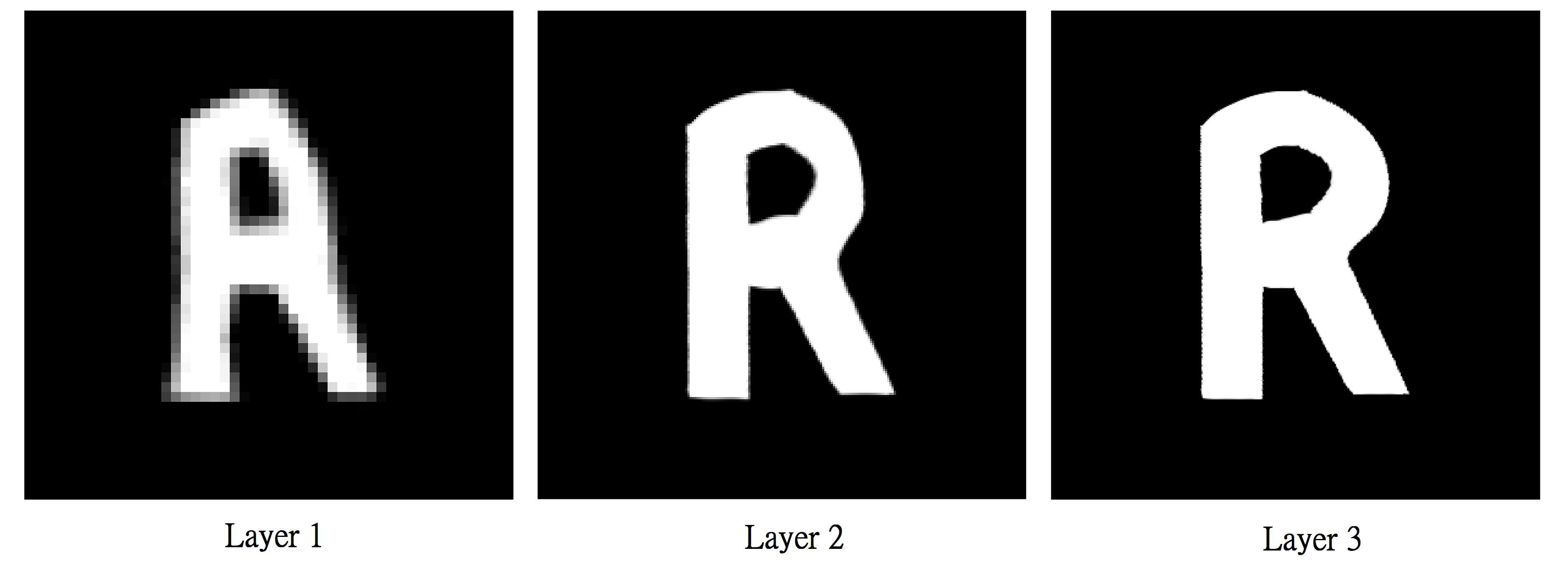}
\caption{The registration results using the multiresolution scheme with 3 layers.}
\end{figure*}

\paragraph{Example 5} We next test our proposed hybrid registration algorithm on a synthetic example. Figure \ref{fig:AR_result}(A) and (B) shows two synthetic images to be registered. (A) shows the image of the character 'A'. (B) shows the image of the character 'R'. Corresponding feature landmarks are labeled on each images. Our goal is to look for a diffeomorphic registration that matches the corresponding landmarks and also the image intensities. (C) shows the obtained diffeomorphic registration using our proposed hybrid registration model. Image 1 in (A) is deformed under the obtain registration to obtain a deformed image, which is shown in (C). The deformed image closely resembles to the target image (Image 2 in (B)). Landmarks are consistently matched, and the obtained registration is bijective. Figure \ref{fig:ARenergy} shows the plots of energy versus iterations. In our algorithm, multi-level approach is applied to perform the registration from the coarsest layer to the finest (original resolution) layer. In this example, three layers are used. Layer 1 refers to the coarsest resolution and layer 3 refers to the finest (original) resolution. The energy plots at different layers are shown in Figure \ref{fig:ARenergy}. The energy is significantly reduced during the optimization process at the first layer. The obtained coarse registration is then interpolate back to the finer layer. An optimal solution is finally reached during the optimization process at the third layer after about 10 iterations. Figure \ref{fig:multiresolutionAR} shows the optimal registration at different layers during the multiresolution scheme.

\begin{figure*}[t]
\centering
\includegraphics[height=1.75in]{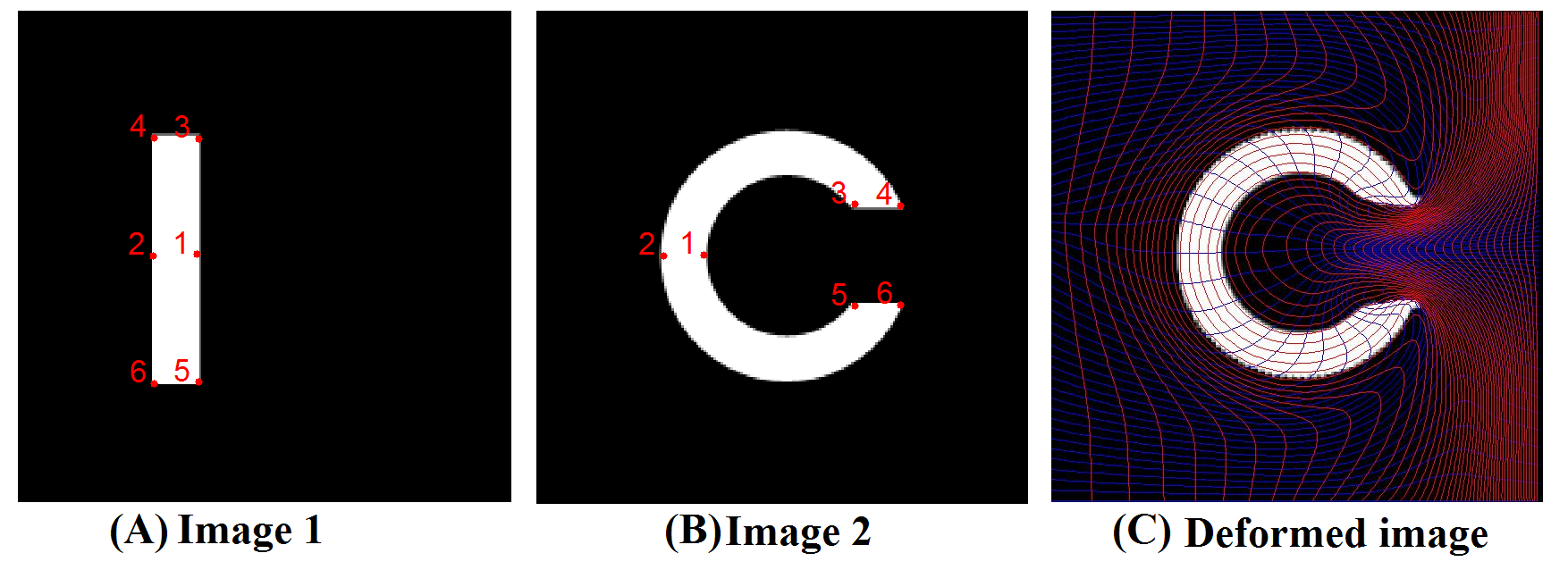}
\caption{(A) shows the image of the character 'I'. (B) shows the image of the character 'C'. These two images are to be registered. Corresponding feature landmarks are labeled on each images. (B) shows the obtained diffeomorphic registration using our proposed hybrid registration model.\label{fig:IC_result}}
\end{figure*}

\paragraph{Example 6} We also test the proposed hybrid registration algorithm on another synthetic example with larger deformation. Figure \ref{fig:IC_result}(A) and (B) shows two synthetic images to be registered. (A) shows the image of the character 'I'. (B) shows the image of the character 'C'. Corresponding feature landmarks are labeled on each images. Again, the goal is to look for a diffeomorphic registration that matches the corresponding landmarks and also the image intensities. (C) shows the obtained diffeomorphic registration using our proposed hybrid registration model. Image 1 in (A) is deformed under the obtain registration to obtain a deformed image, which is shown in (C). The deformed image closely resembles to the target image (Image 2 in (B)). Landmarks are consistently matched, and the obtained registration is bijective.

\begin{figure*}[t]
\centering
\includegraphics[height=2in]{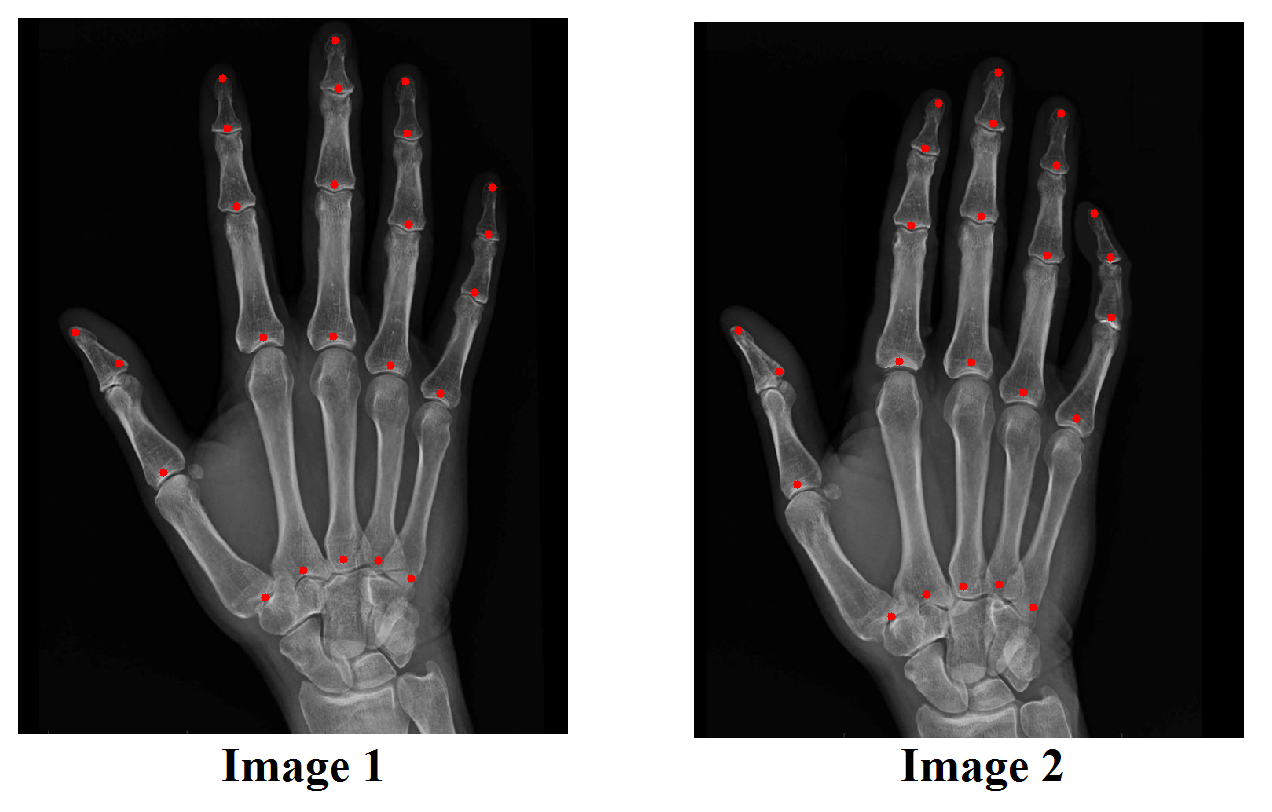}
\caption{Two images of human hands to be registered. Corresponding landmark features are labeled on each images. \label{fig:Hand_register}}
\end{figure*}

\begin{figure*}[t]
\centering
\includegraphics[height=2in]{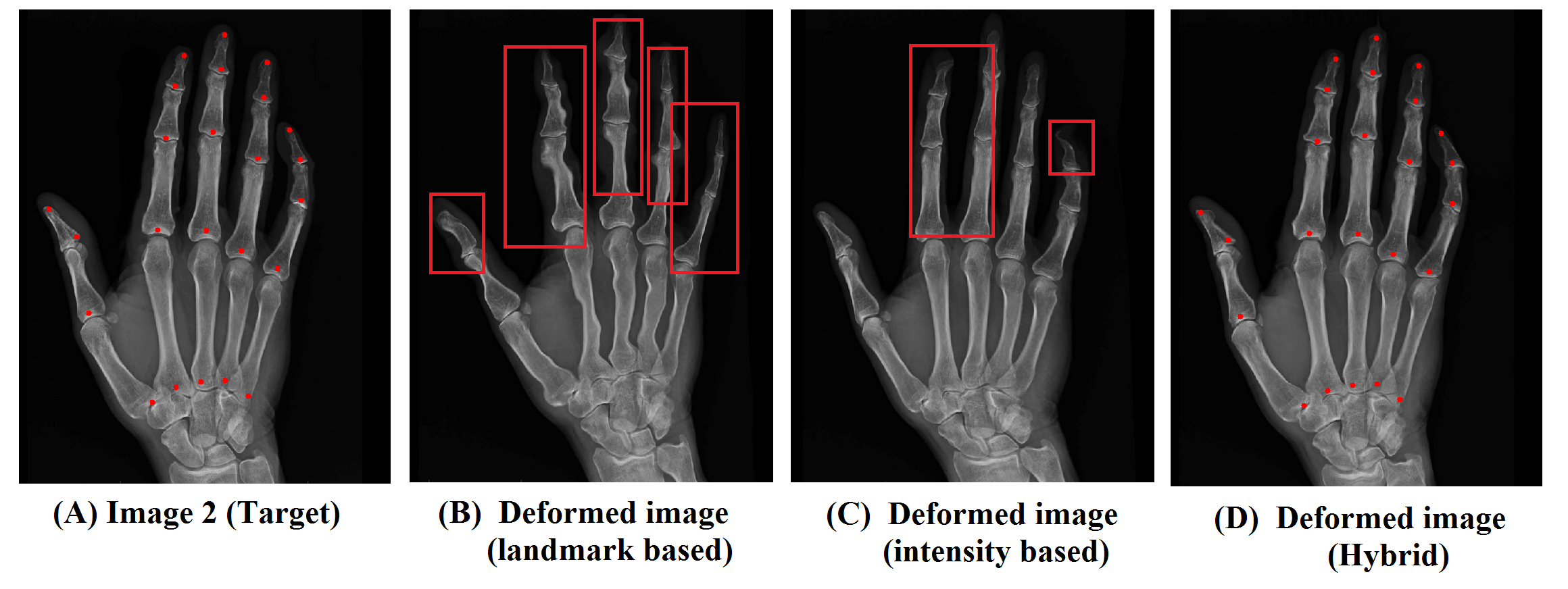}
\caption{Registration results of the human hand images using different approaches. (A) shows the target image (Image 2 as in Figure \ref{fig:Hand_register}). (B) shows the deformed image from Image 1 using the landmark based registration model. (C) shows the deformed image from Image 1 using the intensity based registration model. (D) shows the deformed image from Image 1 using the proposed hybrid registration model. \label{fig:Hand_register_result}}
\end{figure*}

\paragraph{Example 7} We test the hybrid registration algorithm on real images. Figure \ref{fig:Hand_register} shows two images of the human hands. Corresponding landmark features are labeled on each images. In Figure \ref{fig:Hand_register_result}, we show the registration results using different approaches. Figure \ref{fig:Hand_register_result}(B) shows the deformed image from Image 1 using the landmark based registration model. Notice that if we only use landmarks as constraints to guide the registration, the deformed image is very different (see regions in the red boxes) from the target image (as shown in (A)). (C) shows the deformed image from Image 1 using the intensity based registration model. Similarly, the deformed image is very different (see regions in the red boxes) from the target image if only intensity information is used. (D) shows the deformed image from Image 1 using the proposed hybrid registration model. The deformed image closely resembles to the target image. In fact, the intensity mismatching error is less than 1\%, meaning that the registration is very accurate.

\begin{figure*}[t]
\centering
\includegraphics[height=2in]{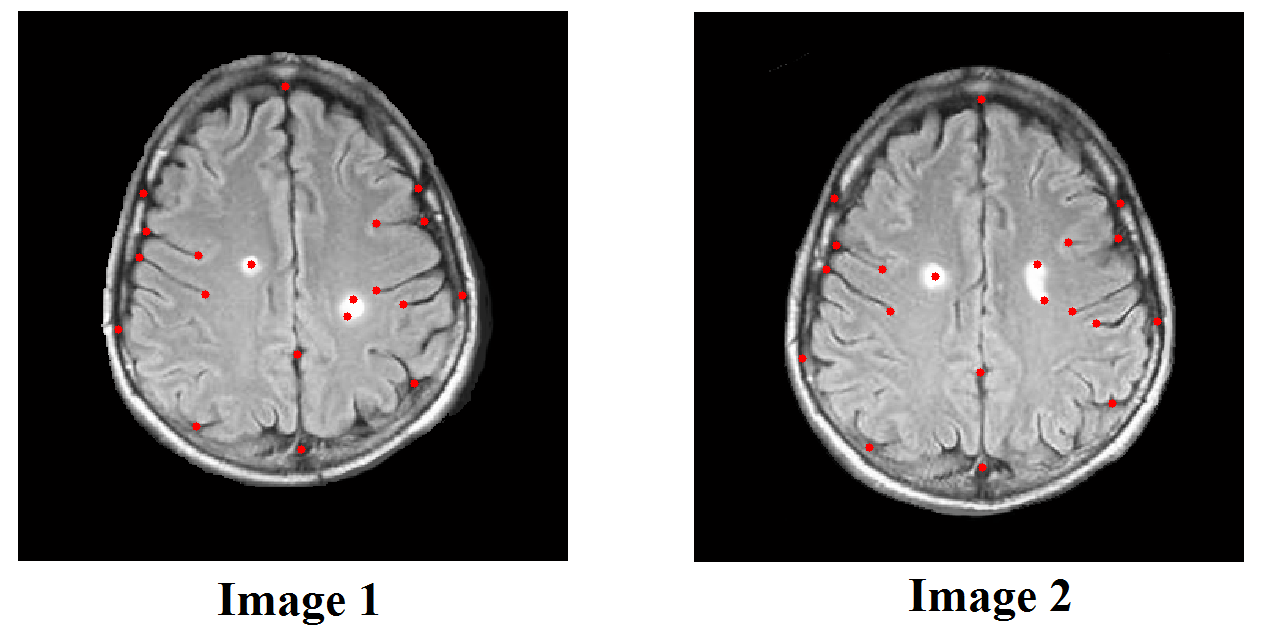}
\caption{Two human brain images to be registered. Corresponding landmark features are labeled on each images. \label{fig:brain_register}}
\end{figure*}

\begin{figure*}[t]
\centering
\includegraphics[height=4.25in]{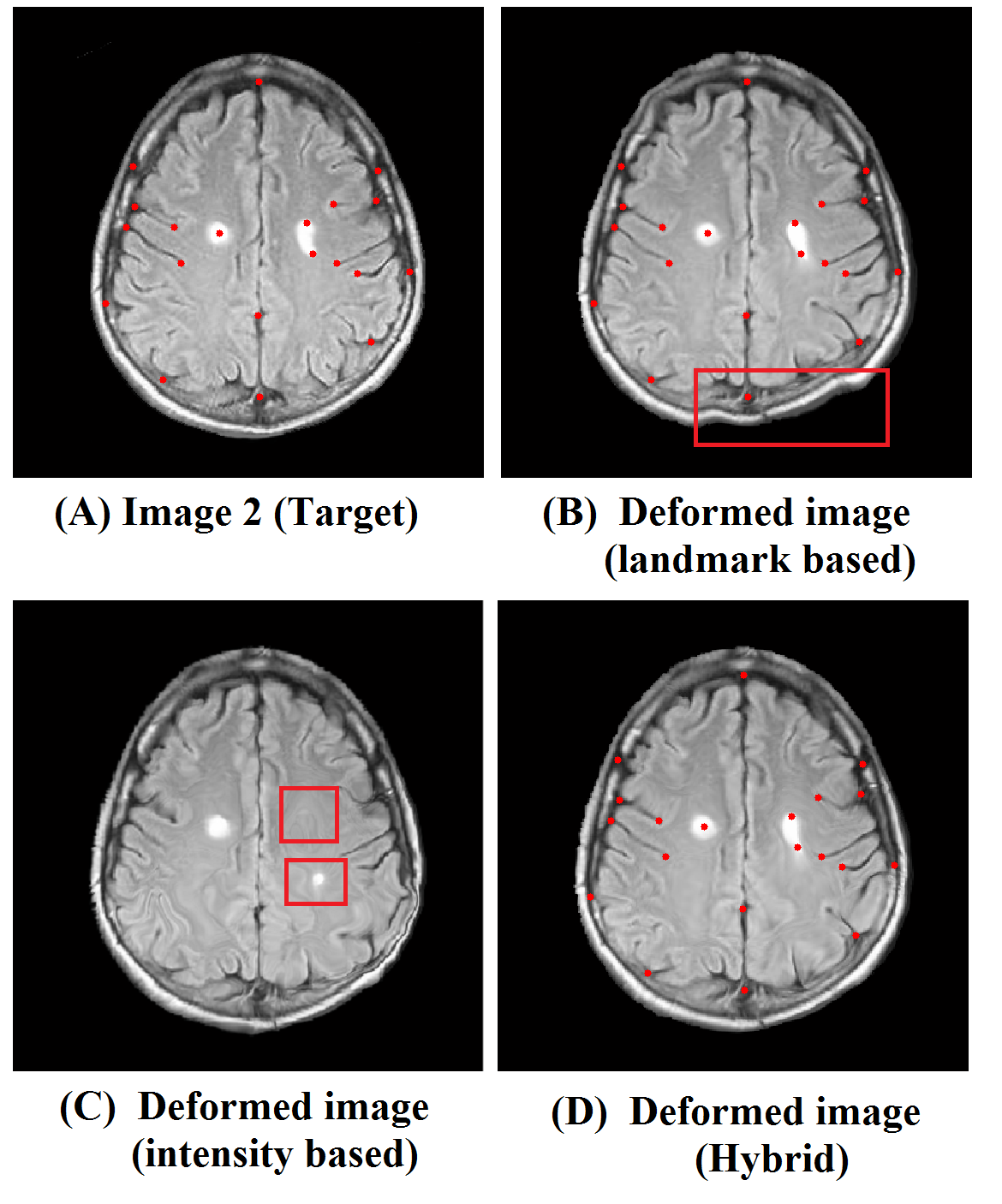}
\caption{Registration results of the human brain images using different approaches. (A) shows the target image (Image 2 as in Figure \ref{fig:brain_register}). (B) shows the deformed image from Image 1 using the landmark based registration model. (C) shows the deformed image from Image 1 using the intensity based registration model. (D) shows the deformed image from Image 1 using the proposed hybrid registration model. \label{fig:brain_registerresult}}
\end{figure*}

\paragraph{Example 8} We also test the hybrid registration algorithm to register two brain MRIs. Figure \ref{fig:brain_register} shows two human brain images. Corresponding landmark features are labeled on each images. In Figure \ref{fig:brain_registerresult}, we show the registration results using different approaches. Figure \ref{fig:brain_registerresult}(B) shows the deformed image from Image 1 using the landmark based registration model. If we only use landmarks as constraints to guide the registration, the deformed image is different (see regions in the red boxes) from the target image (as shown in (A)). (C) shows the deformed image from Image 1 using the intensity based registration model. Similarly, the deformed image is very different (see regions in the red boxes) from the target image if only intensity information is used. Notice that the tumor in the red box does not move using intensity based registration. (D) shows the deformed image from Image 1 using the proposed hybrid registration model. The deformed image closely resembles to the target image. The intensity mismatching error is less than 1\%, meaning that the registration is very accurate.

\begin{figure*}[t]
\centering
\includegraphics[height=2in]{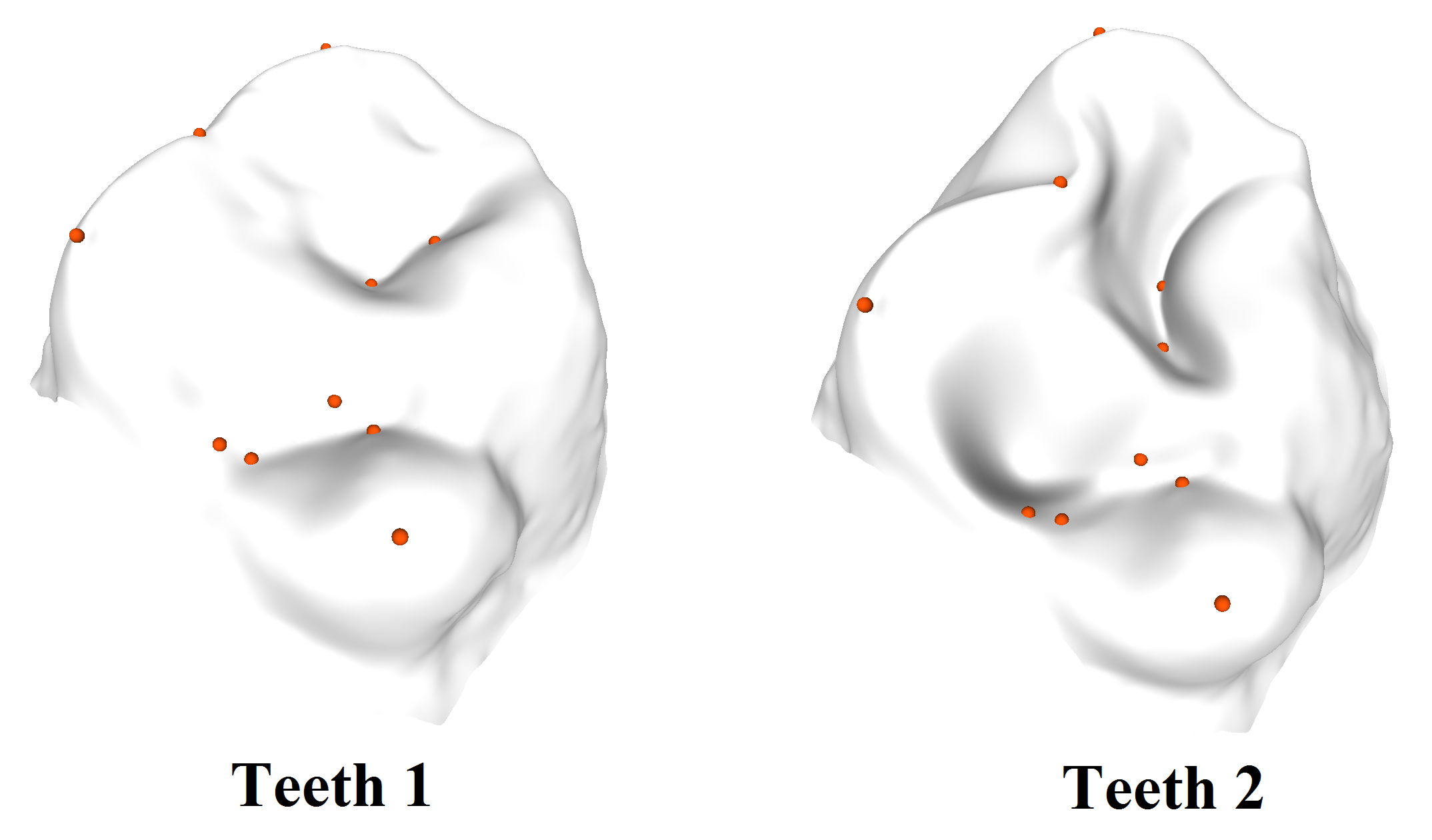}
\caption{Two human teeth to be registered, each of them are labeled with corresponding landmarks.\label{fig:Teeth12}}
\end{figure*}

\begin{figure*}[t]
\centering
\includegraphics[height=1.5in]{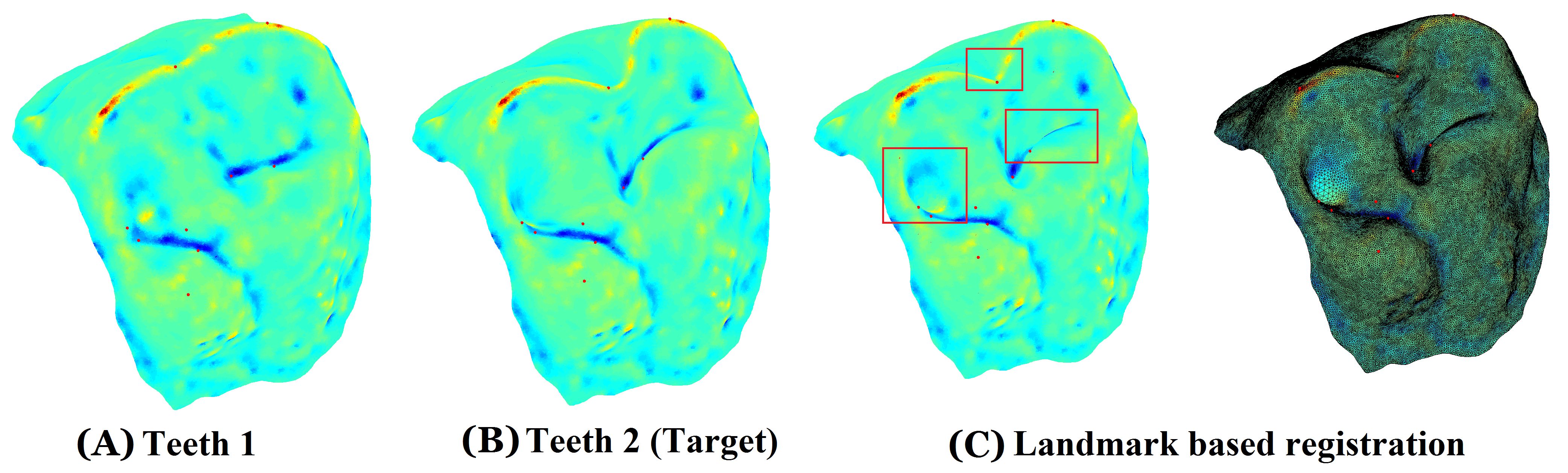}
\caption{Registration results of the two teeth surfaces using the landmark-matching quasi-conformal registration. (A) shows the surface of Teeth 1, whose colormap is given by its mean curvature. (B) shows the surface of Teeth 2, whose colormap is given by its mean curvature. (C) shows the registration result using the landmark-matching quasi-conformal registration. The colormap on the surface of Teeth 1 is mapped to the surface of Teeth 2 using the obtained registration. Note that the curvature is not matched consistently (see the regions in the red boxes).  \label{fig:Teeth_registration_lm}}
\end{figure*}

\begin{figure*}[t]
\centering
\includegraphics[height=1.5in]{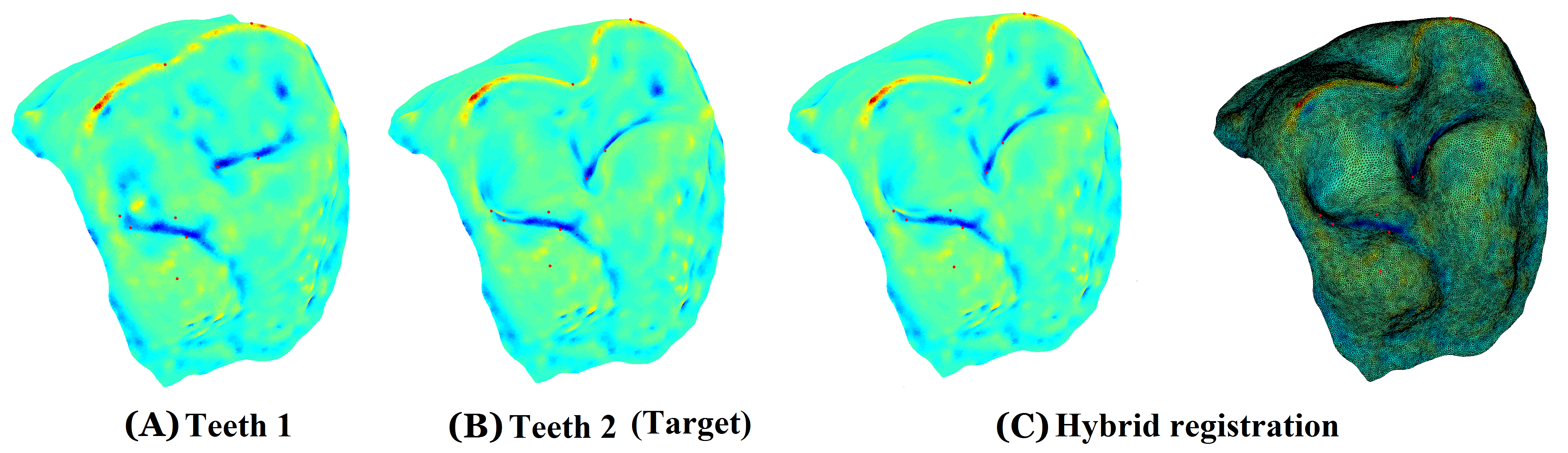}
\caption{Registration results of the two teeth surfaces using the hybrid quasi-conformal registration. (A) shows the surface of Teeth 1, whose colormap is given by its mean curvature. (B) shows the surface of Teeth 2, whose colormap is given by its mean curvature. (C) shows the registration result using the hybrid quasi-conformal registration. The colormap on the surface of Teeth 1 is mapped to the surface of Teeth 2 using the obtained registration. Note that the curvature is matched consistently.  \label{fig:Teeth_registration_summary}}
\end{figure*}

\paragraph{Example 9} We also test the hybrid registration algorithm to register two human teeth surfaces. Figure \ref{fig:Teeth12} shows two human teeth surfaces, each of them are labeled with corresponding landmarks. Figure \ref{fig:Teeth_registration_lm} shows the registration results of the two teeth surfaces using the landmark-matching quasi-conformal registration. (A) shows the surface of Teeth 1, whose colormap is given by its mean curvature. (B) shows the surface of Teeth 2, whose colormap is given by its mean curvature. (C) shows the registration result using the landmark-matching quasi-conformal registration. The colormap on the surface of Teeth 1 is mapped to the surface of Teeth 2 using the obtained registration. Note that the curvature is not matched consistently (see the regions in the red boxes). It shows that the registration is not accurate if only landmarks are used to guide the registration process. Figure \ref{fig:Teeth_registration_summary} (C) shows the registration result using the hybrid quasi-conformal registration. In this case, both landmarks and curvature information are used to guide the registration process. The colormap (mean curvature) on the surface of Teeth 1 is mapped to the surface of Teeth 2 using the obtained registration. Note that the curvature is matched consistently, which means the registration result is accurate.

\begin{figure*}[t]
\centering
\includegraphics[height=2in]{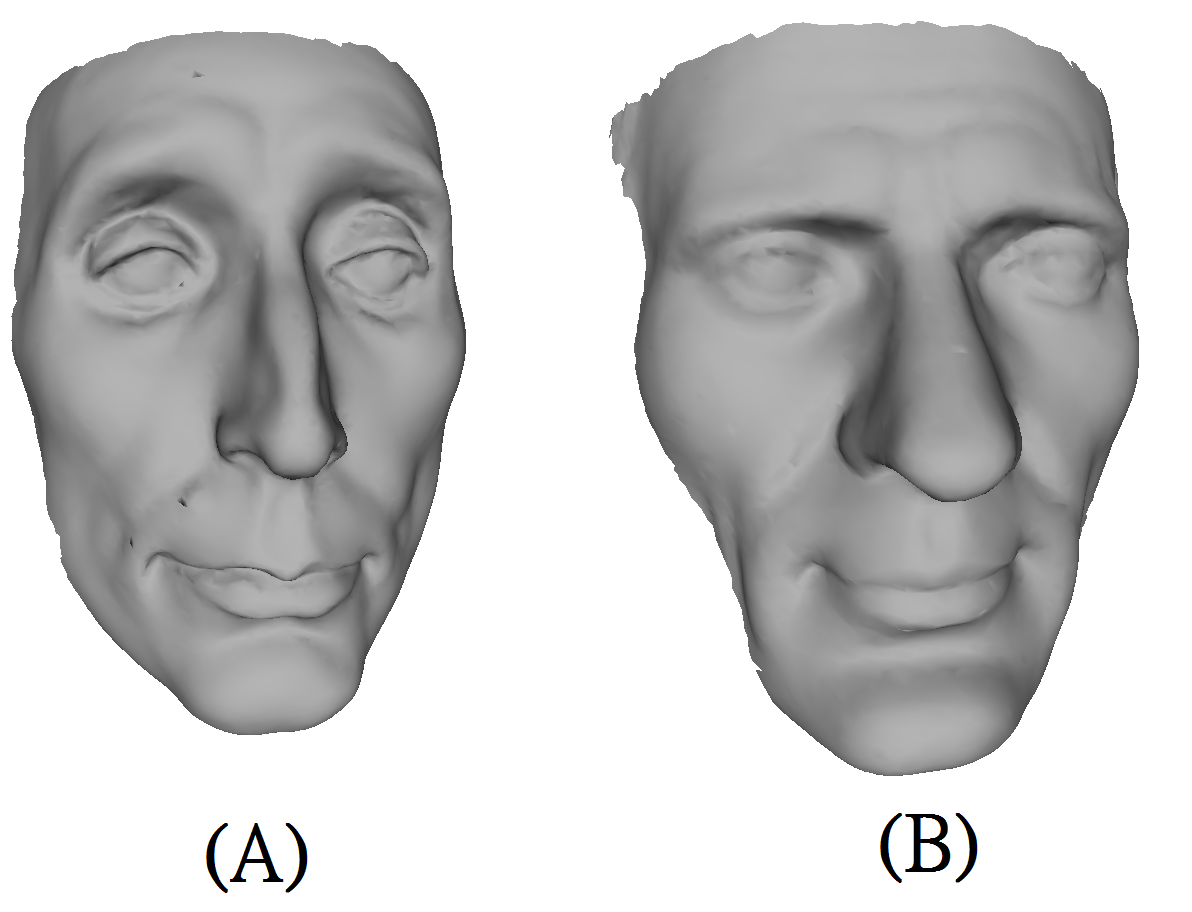}
\caption{Two human faces to be registered, each of them are labeled with corresponding landmarks. The mesh data are freely available at  http://shapes.aimatshape.net. \label{fig:Face_mesh}}
\end{figure*}

\begin{figure*}[t]
\centering
\includegraphics[height=2in]{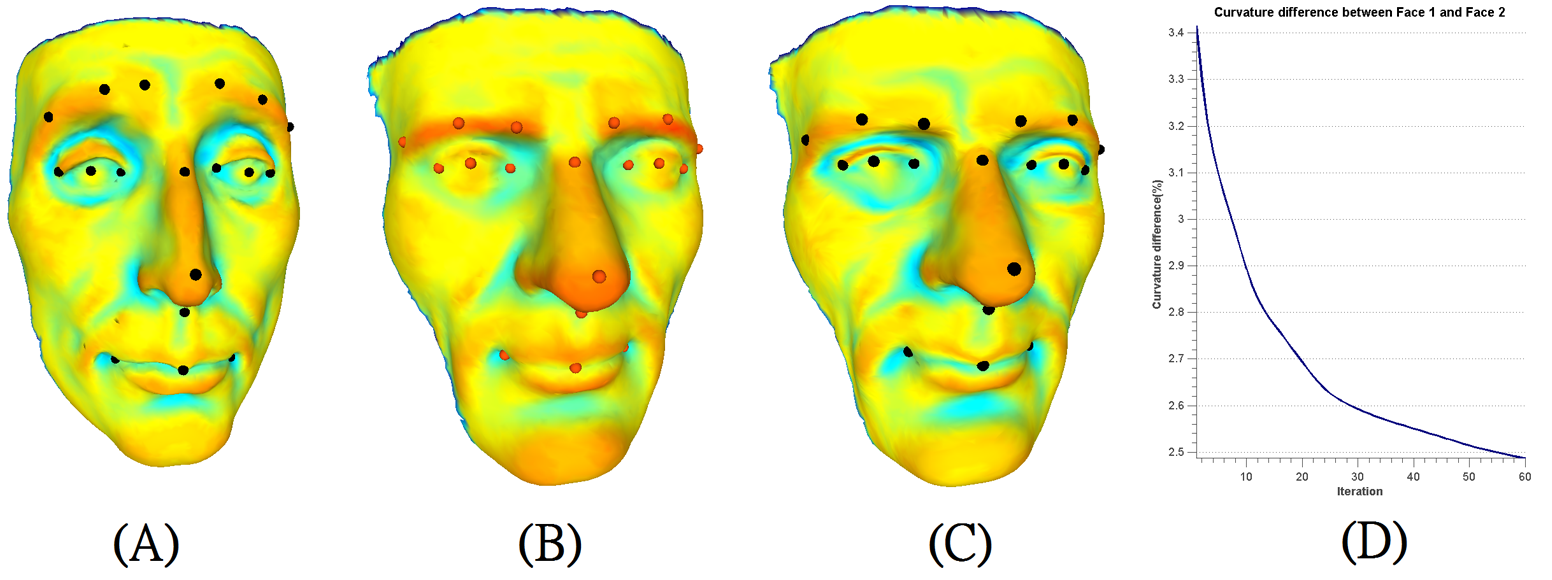}
\caption{Registration results of the two human faces using the hybrid quasi-conformal registration. (A) shows the surface of human face 1, whose colormap is given by its mean curvature. (B) shows the surface of human face 2, whose colormap is given by its mean curvature. (C) shows the registration result using the hybrid quasi-conformal registration. The colormap on the surface of human face 1 is mapped to the surface of human face 2 using the obtained registration. Note that the corresponding regions are consistently matched. (D) shows the plot of curvature mismatching versus iterations.\label{fig:Face_matching}}
\end{figure*}

\paragraph{Example 10} We also test the hybrid registration algorithm to register two human face surfaces. Figure \ref{fig:Face_mesh} shows two human face surfaces, each of them are labeled with corresponding landmarks. Figure \ref{fig:Face_matching} shows the registration results of the two human faces using the landmark-matching quasi-conformal registration. (A) shows the surface of human face 1, whose colormap is given by its mean curvature. (B) shows the surface of human face 2, whose colormap is given by its mean curvature. (C)  shows the registration result using the hybrid quasi-conformal registration. The colormap on the surface of human face 1 is mapped to the surface of human face 2 using the obtained registration. Note that the corresponding regions are consistently matched. (D) shows the plot of curvature mismatching versus iterations. It shows that our algorithm iteratively adjust the quasi-conformal registration to an optimal one that minimizes the curvature mismatching error.

\section{Conclusion}\label{conclusion}
This paper presents a novel method to obtain diffeomorphic image or surface registrations with large deformations via quasi-conformal maps. The main strategy is to minimize an energy functional involving a Beltrami coefficient
term. The Beltrami coefficient measures the conformality distortion of the quasi-conformal map. It controls controls the bijectivity and smoothness of the registration. By minimizing the energy functional, we obtain an optimal
Beltrami coefficient associated to the desired registration, which is guaranteed to be bijective, even with very large deformations. The proposed method can be applied for both landmark based registration and hybrid registration. Experiments have been carried out on both synthetic and real data. Results show that our proposed method can effectively obtain diffeomorphic registration between images or surfaces with least amount of local geometric distortion. The obtained registration is guaranteed to be bijective (1-1 and onto), even with a large deformation or large number of landmark constraints. In the future, we plan to extend the proposed method to high-genus surfaces and apply the method to more real applications in medical imaging for disease analysis.

\end{document}